\documentclass[10pt]{amsart} 
\usepackage{amsmath,amssymb}
\usepackage{amsrefs}

\numberwithin{equation}{section}
\newtheorem{theorem}[equation]{Theorem}
\newtheorem{lemma}[equation]{Lemma}
\newtheorem{corollary}[equation]{Corollary}
\newtheorem{prop}[equation]{Proposition}

\theoremstyle{definition}
\newtheorem{definition}[equation]{Definition}
\newtheorem{notation}[equation]{Notation}
\newtheorem{example}[equation]{Example}

\theoremstyle{remark}
\newtheorem{remark}[equation]{Remark}

\begin{document}

\begin{center}
\texttt{Comments, suggestions, corrections, and references 
  welcomed!}
\end{center}

\title[Capable groups of prime exponent and class two~II]{Capable
  groups of prime exponent and class two~II}
\author{Arturo Magidin} 
\address{Dept.~Mathematical Sciences, The University of Montana, Missoula MT 59812-0864}
\curraddr{Mathematics Dept. University of Louisiana at Lafayette,
  217 Maxim Doucet Hall, P.O.~Box 41010, Lafayette LA 70504-1010}
\email{magidin@member.ams.org}

\subjclass[2000]{Primary 20D15, Secondary 20F12, 15A04}

\begin{abstract}
We consider the capability of $p$-groups of class two and
odd prime exponent. We use linear algebra and counting arguments
to establish a number of new results. In particular, we
settle the $4$-generator case, and prove a sufficient condition
based on the ranks of $G/Z(G)$ and $[G,G]$.
\end{abstract}

\maketitle

\section*{Introduction.}

In this work we continue to study the capability of finite
$p$-groups of class two and prime exponent, using the approach 
introduced in~\cite{capablep}. We restrict to odd primes,
since the case of abelian groups is well understood. I owe a 
considerable amount of the material here to discussions with David
McKinnon\nocite{persdave}, who helped me clarify the ideas and
helped with most of the geometry.

A group $G$ is said to be \textit{capable} if and only if there exists
a group $K$ such that $G=K/Z(K)$. For $p$-groups (groups of finite
prime-power order) capability is closely related to their
classification. Baer characterized the capable groups which are direct
sums of cyclic groups in~\cite{baer}; the capable extra-special
$p$-groups were characterized by Beyl, Felgner, and Schmid
in~\cite{beyl} (only the dihedral group of order~$8$ and the
extra-special groups of order $p^3$ and exponent~$p$ are capable);
they  also described the metacyclic groups which are capable. The
author characterized the $2$-generated capable $p$-groups of class
two~\cites{capable,twocubed} (for odd~$p$, independently obtained in
part by Bacon and~Kappe in~\cite{baconkappe}).

For the case of $p$-groups of class two and exponent $p$, where $p$ is
an odd prime, some necessary and some sufficient conditions for
capability are known, and so there is a hope expressed
in~\cite{baconkappe} that a full characterization for this class may
be tractable with current techniques. We began to study this situation
in~\cite{capablep}; there we described a way to translate the problem
into a statement of linear algebra, and obtained several results using
this restatement. We will introduce what I hope is clearer notation
here and obtain new results, including a new sufficient condition. We
direct the reader to~\cite{capablep} for basic definitions and
notation.

In Section~\ref{sec:setup} we recap the characterization of
capability first described in~\cite{capablep} (along the way we will
correct a slight misstatement made in the proof there), and then describe
the linear algebra problem that is derived from it.
In Section~\ref{sec:linearalgebra}
we prove some basic facts about the linear algebra situation. In
Section~\ref{sec:dimcounting} we use a counting argument to derive a
sufficient condition for the capability of~$G$, based on the ranks of
$G^{\rm ab}$ and $[G,G]$. Then in Section~\ref{sec:oldresults} we recall
some results proven in~\cite{capablep} and obtain some new ones. We
will note there the relevance of element of $G$ which lie in $Z(G)$ and
have nontrivial image in $G^{\rm ab}$. In
Section~\ref{sec:centralelements} we will therefore establish a
connection between these elements and certain subspaces discussed in
Section~\ref{sec:dimcounting}.
Finally, in
Section~\ref{sec:nequalfour} we settle the four-generated case for
this class.

\section{The set-up.}\label{sec:setup}

Throughout $p$ will be used to denote an odd prime.
We begin by describing a ``canonical witness'' to the capability of a
finite group $G$ of class two and odd prime exponent. For the
definition of the $3$-nilpotent product, see~\cites{capable,struikone}.

\begin{theorem}[Theorem~2.3 in~\cite{capablep}]
Let $G$ be a finite noncyclic group of class at most two and exponent an odd
prime~$p$. Let $g_1,\ldots,g_n$ be elements of $G$ that project onto a
basis for $G^{\rm ab}$, and let $F$ be the $3$-nilpotent product of
$n$ cyclic groups of order~$p$, generated by $x_1,\ldots,x_n$
respectively. Let $N$ be the kernel of the morphism $F\to G$
induced by mapping $x_i\mapsto g_i$ for each $i$. Then $G$ is capable
if and only if
\[ G \cong \left(F/[N,F]\right)\bigm/ Z\left(F/[N,F]\right).\]
\label{th:basicequivalence}
\end{theorem}

\begin{proof}
Note that since $[N,F]\subset F_3$ and $N[N,F]\subset Z(F/[N,F])$, it is
always the case that the central quotient of $F/[N,F]$ is a quotient
of~$G$.

The ``if'' clause of the theorem immediate. For the converse, assume
that $G$ is capable, and let $K$ be a group with $K/Z(K)\cong G$. Let
$k_1,\ldots,k_n$ be elements of $K$ that project onto
$g_1,\ldots,g_n$, respectively. The subgroup of~$K$ generated by
$k_1,\ldots,k_n$ is cocentral, hence has central quotient isomorphic
to $K/Z(K)$; we may thus assume that $K$ is generated by
$k_1,\ldots,k_n$. 

Note that $K$ is of class at most $3$. We claim that
$K_2^p=\{e\}$. Indeed, since $K/Z(K)$ is of exponent $p$, it follows
that the $p$-th power of any element of $K$ is in $Z(K)$. So, if $c\in
K_2$ and $k\in K$, then 
\[ e = [c,k^p] = [c,k]^p [c,k,k]^{\binom{p}{2}} = [c,k]^p.\]
In particular, $K_3$ is of exponent~$p$. If $h,k\in K$, then we
have:
\[ e = [h,k^p] = [h,k]^p[h,k,k]^{\binom{p}{2}} = [h,k]^p,\]
so every commutator is of exponent~$p$. Since $K_2$ is abelian, it
follows that $K_2$ is of exponent~$p$, as claimed.

Let $\mathfrak{F}$ be the $3$-nilpotent product of $n$ infinite cyclic
groups, and let $\mathcal{F}=\mathfrak{F}/(\mathfrak{F}_2)^p$. Denote
the generators of the cyclic groups by $y_1,\ldots,y_n$. Then
$\mathcal{F}$ is the relatively free group of rank $n$ in the variety
of all groups of class at most $3$ in which the commutator subgroup is
of exponent~$p$. Note that $F\cong \mathcal{F}/\langle
y_1^p,\ldots,y_n^p\rangle$, and that under this isomorphism we have an
identification of $F_2$ with $\mathcal{F}_2$. Let $\mathcal{N}$ be the
subgroup of $\mathcal{F}_2$ corresponding to the subgroup $N$ of $F$.

The map sending $y_i$ to $k_i$ induces a morphism $\mathcal{F}\to
K$. Since $\mathcal{N}$ maps to $e$ under the composite map
$\mathcal{F}\to K \to K/Z(K)$, it follows that $\mathcal{N}$ maps into
$Z(K)$, and hence that the map $\mathcal{F}\to K$ factors through
$\mathcal{F}/[\mathcal{N},\mathcal{F}]$. As this map sends the center
of the latter group into the center of $K$, it follows that $G$ is a
quotient of $(\mathcal{F}/[\mathcal{N},\mathcal{F}]) /
Z(\mathcal{F}/[\mathcal{N},\mathcal{F}])$. 

Using the normal forms for $\mathcal{F}$ (see \cite{struikone}), it is
easy to show that the central quotients of
$\mathcal{F}/[\mathcal{N},\mathcal{F}]$ and of $\mathcal{F}/\langle
y_1^p,\ldots,y_n^p,[\mathcal{N},\mathcal{F}]\rangle$ are
isomorphic. This latter group is of course isomorphic to $F/[N,F]$, so
we conclude that if $G$ is capable, then it is a quotient of
$(F/[N,F])/Z(F/[N,F])$.

So if $G$ is capable, then $G$ is isomorphic to a quotient of the
central quotient of $F/[N,F]$ (which is a finite group), and in turn
has the central quotient of $F/[N,F]$ as a quotient. The only way this
is possible is for $G$ to be isomorphic to the central quotient of
$F/[N,F]$, as claimed.
\end{proof}

\begin{remark} In \cite{capablep} there was a slight error in the
  assertion just prior to Theorem~2.3. It was asserted there that 
  $Z(F/(F^p[N,F])) \cong Z(F/[N,F])/F^p$, but this assertion is false
  when $p=3$ (it is correct if $p>3$ by regularity). The
  subgroup $F^p$ should be replaced by the subgroup generated by the
  $p$-th powers of the generators, as we did in the proof above.
\end{remark}

Since $G$ is isomorphic to $F/NF_3$ and $F_3\subset Z(F)$, the
following result now follows:

\begin{theorem}[Theorem~3.6 in~\cite{capablep}] 
Let $G$ be a finite group of class at most two and exponent an odd
prime~$p$. Let $g_1,\ldots,g_n$ be elements of $G$ that project onto a
basis for $G^{\rm ab}$, and let $F$ be the $3$-nilpotent product of
$n$ cyclic groups of order~$p$, generated by $x_1,\ldots,x_n$
respectively. Let $N$ be the kernel of the canonical map $F\to G$
induced by mapping $x_i\mapsto g_i$ for each $i$, and let $C$ be the
subgroup of $F_2$ spanned by all basic commutators of the form
$[x_j,x_i]$, $1\leq i<j\leq n$. Write $N=X\oplus
F_3$, where $X\subset C$.
Then $G$ is capable
if and only if
$\bigl\{ x \in C \,\bigm|\,
   [x,F]\subset [X,F]\bigr\} = X$.
\label{th:basicequivalencewithN}
\end{theorem}

It is at this point that we introduce linear algebra as a way to
codify the situation we wish to study: we have
abelian group homomorphisms between elementary abelian $p$-groups
given by $[-,x_k]\colon F_2\to F_3$. Thus, we may
interpret it as linear maps between vectors spaces over
$\mathbb{F}_p$. If $X$ is a subgroup of $C$, then $[X,F]$ will be
the span of the images of $[X,x_k]$ with $k=1,\ldots,n$. On the other
hand, the elements $x\in C$ such that
$[x,F]$ is contained in~$[X,F]$ correspond to those elements whose
images under all maps $[-,x_k]$ lie in $[X,F]$, i.e., the intersection
of the pullbacks of these maps. This suggests
the following general construction and definitions:

\begin{definition} Let $V_1$ and $V_2$ be vector spaces,
  and let $\{\ell_i\colon V_1\to V_2\}_{i \in I}$ be a nonempty family of
  linear transformations. 
  For every subspace $X$ of $V_1$, we define ${X}^* < V_2$ to be the
  subspace spanned by all the images of $X$ under the linear
  transformations $\ell_i$, that is:
\[ X^* = \Bigl\langle \ell_i(X) \Bigm| i\in I\Bigr\rangle.\]
For every subspace $Y$ of $V_2$, we define $Y^* < V_1$ to be the
intersection of all pullbacks of $Y$ under the transformations
$\ell_i$, that is:
\[ Y^* = \bigcap_{i\in I} \ell_i^{-1}(Y).\]
\end{definition}

Note that for all subspaces of $V_1$, if $X\subset X'$, then
$X^*\subset X'^*$; likewise, for subspaces of $V_2$, if $Y\subset
Y'$, then $Y^*\subset Y'^*$.

\begin{theorem} Let $V_1$ and $V_2$ be vector spaces, 
  and let $\{\ell_i\colon V_1\to
  V_2\}_{i\in I}$ be a nonempty family of linear transformations. The operator on the
  subspaces of $V_1$ defined by $X\mapsto X^{**}$ is a closure
  operator; that is, it is increasing, isotone, and idempotent.
  Moreover, if
  $X$ is a subspace of $V_1$, then $X^* = X^{***}$.
\end{theorem}

\begin{proof} It is clear that
  $X\subset X^{**}$, so the operator is increasing; since $X\subset
  X'$ implies $X^*\subset X'^*$, which in turn implies $X^{**}\subset
  X'^{**}$, the operator is isotone.

  To show the operator is idempotent, we want to show that
  $X^{**} = (X^{**})^{**}$. For simplicity, write $X^{**}=Z$. Since
  the operator is increasing, we know that $Z\subset Z^{**}$. By
  construction, we also have that $\ell_i(Z)\subset X^*$ for each $i$,
  so $Z^*\subset X^*$. From this, we get that $Z^{**} \subset
  X^{**}=Z$, as desired.

  Finally, let $X < V_1$. Since $X\subset X^{**}$, we must have $X^*\subset
  X^{***}$. Conversely, from the definition of $X^{***}$ it follows
  that $\ell_i(X^{**})\subset X^*$ for all $i$, so $X^{***}\subset X^*$,
  giving equality.
\end{proof}

\begin{remark}
It may be worth noting that this closure operator is
algebraic (meaning that the closure of a subspace $X$ is the
union of the closure of all finitely generated subspaces $X'$
contained in~$X$), though it is not topological (in general,
the closure of the subspace generated by $X$ and $Y$ is not the
subspace generated by $X^{**}$ and~$Y^{**}$). 

\end{remark}

The dual property is true for subspaces of $V_2$:

\begin{theorem} Let $V_1$ and $V_2$ be vector spaces, 
  and let $\{\ell_i\colon V_1\to V_2\}_{i\in I}$ be a
  nonempty family of
  linear transformations.
 The operator on subspaces of $V_2$ defined by $Y\mapsto
  Y^{**}$ is an interior operator; that is, it is decreasing, isotone,
  and idempotent. Moreover, if $Y$ is a subspace of $V_2$, then $Y^*=Y^{***}$.
\end{theorem}

\begin{proof} As before, the operator is isotone. Since
  $Y^{**} = \langle \ell_i(Y^*)\rangle$, and 
  $Y^*\subset \ell_i^{-1}(Y)$ for each $i$, it follows that
  $Y^{**}\subset Y$, so the operator is decreasing.

To show the operator is idempotent, set $Z=Y^{**}$. Since $Y^*\subset
Y^{***}$, we must have $Z = Y^{**}\subset Y^{****}=Z^{**}$.
The
reverse inclusion always holds, so $Z=Z^{**}$, as desired.

Finally, since $Y^{**}\subset Y$, it follows that $Y^{***}\subset
Y^*$.  But since the operator ${}^{**}$ on subspaces of $V_1$ is
increasing, we also know that $Y^*\subset Y^{***}$, giving equality.
\end{proof}

\begin{definition} Let $V_1$, $V_2$ be vector spaces, and let
  $\{\ell_i\colon V_1\to V_2\}_{i\in I}$ be a family of linear
  transformations. We will say that a subspace $X$ of $V_1$ is
  \textit{$\{\varphi_i\}_{i\in I}$-closed} (or simply \textit{closed} if there is no
  danger of ambiguity) if and only
  if $X=X^{**}$. 
\end{definition}

To tie in this construction with our capability problem, we define
specific spaces and maps. We fix an odd prime $p$ throughout. Let $n$
be an integer, $n>1$.  We define three vectors spaces with
distinguished bases:

\begin{definition}
Let $U(n)$ be the vector space over the finite field $\mathbb{F}_p$ of
$p$ elements, with basis vectors $u_1,\ldots,u_n$.  Let $V(n)$ be the
vector space over $\mathbb{F}_p$ of dimension $\binom{n}{2}$, with
basis vectors $v_{ji}$, $1\leq i<j\leq n$. Let $W(n)$ be the vector
space over $\mathbb{F}_p$ of dimension $2\binom{n}{2} + 2\binom{n}{3}
= 2\binom{n+1}{3}$, with basis vectors $w_{jik}$, $1\leq i<j\leq n$,
$i\leq k\leq n$.
\end{definition}

\begin{notation} For simplicity, we will sometimes use $v_{ij}$ with
  $i< j$ for an element of $V(n)$; in that case, we understand it
  to mean that  $v_{ij}=-v_{ji}$. Likewise, $v_{ii}=\mathbf{0}$ for
  each~$i$.
\end{notation}

If there is no danger of ambiguity, and $n$ is understood from
context, we will simply use $U$, $V$, and~$W$ instead of $U(n)$,
$V(n)$, and $W(n)$. This will be the case in most of our applications.

We define two families of $n$ linear operators:

\begin{definition}
The linear maps $\psi_i\colon U\to V$, $i=1,\ldots,n$ are defined to be:
\[ \psi_i(u_j) = v_{ji} = \left\{\begin{array}{ll}
v_{ji}&\mbox{if $j>i$}\\
-v_{ij}&\mbox{if $j<i$}\\
\mathbf{0}&\mbox{if $i=j$.}
\end{array}\right.\]
\end{definition}

\begin{definition}
The linear maps $\varphi_k\colon V \to W$, $k=1,\ldots,n$,
are defined to be:
\[ \varphi_k(v_{ji}) = \left\{\begin{array}{ll}
w_{jik} & \mbox{if $k\geq i$,}\\
w_{jki}-w_{ikj} & \mbox{if $k<i$.}
\end{array}\right.\]
\end{definition}

Note that there is a natural identification of $V$ with $U\wedge U$,
and of $W$ with $(U\wedge U) \otimes U$ modulo the Jacobi identity:
\[ (u_i\wedge u_j) \otimes u_k + (u_j\wedge u_k)\otimes u_i +
(u_k\wedge u_i)\otimes u_j = 0.\]
Under these identifications, the maps $\psi_i$ correspond to
$\psi_i(\mathbf{u}) = \mathbf{u}\wedge u_i$, and the maps $\varphi_k$
correspond to $\varphi_k(\mathbf{v}\wedge \mathbf{w}) =
(\mathbf{v}\wedge \mathbf{w})\otimes u_k$. 

The connection between capability and this construction is now
apparent: $U$ corresponds to $F^{\rm ab}$ by identifying the
generators $u_i$ with $x_i$, $i=1,\ldots,n$. The space $V$ corresponds
to the subgroup $\langle [x_j,x_i]\rangle$ by identifying $v_{ji}$
with the commutator $[x_j,x_i]$, $1\leq i<j\leq n$. And $W$ corresponds to $F_3$ by
identifying $w_{jik}$ with $[x_j,x_i,x_k]$, $1\leq i<j\leq n$, $i\leq
k\leq n$. The maps $\psi_i$ correspond to maps $[-,x_i]\colon F/F_2
\to F_2/F_3$, while $\varphi_k$ corresponds to the map $[-,x_k]\colon
F_2/F_3 \to F_3$. In particular, if $X$ is a subspace of $V$, then it
corresponds to a subgroup $N$ of $F$ contained in $\langle
[x_j,x_i]\rangle$; and by Theorem~\ref{th:basicequivalencewithN}, it
follows that the group $F/NF_3$ is capable if and only if
$X=X^{**}$. Therefore, we have:

\begin{theorem}
Let $G$ be a finite noncyclic group of class at most two and exponent an odd
prime~$p$. Let $g_1,\ldots,g_n$ be elements of $G$ that project onto a
basis for $G^{\rm ab}$, and let $F$ be the $3$-nilpotent product of
$n$ cyclic groups of order~$p$, generated by $x_1,\ldots,x_n$
respectively. Let $N$ be the kernel of the morphism $F\to G$
induced by mapping $x_i\mapsto g_i$ for each $i$. Let $N=X\oplus
F_3$, where $X\subset \langle [x_j,x_i]\,|\, 1\leq i<j\leq n\rangle$,
and identify $X$ with the corresponding subspace of $V(n)$. 
Then $G$ is capable
if and only if $X$ is closed with respect to $\{\varphi_i\}_{i=1}^{n}$.
\label{th:equivwithclosure}
\end{theorem}

Thus, the question of which $n$-generated $p$-groups of class at most
$2$ and exponent~$p$ ($p$ odd) are capable is equivalent to the
question of which subspaces of $V(n)$ are closed. Note, however, that
distinct subspaces may correspond to isomorphic groups: for example,
if we permute the indices, we obtain a different $X$ but an
isomorphic group~$G$.

The reason we include the maps $\psi_i$ may seem a bit more
mysterious. Note that if $Z$ is a subspace of $U$, then identifying
$V$ with $U\wedge U$ we have that $Z^*$ is none other than $Z\wedge
U$. This corresponds to the relations necessary to make the elements
of~$F$ corresponding to $Z$ central in~$F$; these relations have an
important interplay with the morphisms $\varphi_i$ and with questions
of capability, as we will see in Sections~\ref{sec:oldresults}
and~\ref{sec:centralelements}.

\section{The linear algebra.}\label{sec:linearalgebra}

We collect here some observations on the spaces $V$ and~$W$, and the
linear transformations $\varphi_k$ defined above.

\begin{definition} Let $i,j$ be integers, $1\leq i< j \leq n$. We
  let $\pi_{ji}\colon V\to \langle v_{ji}\rangle$ be the
  canonical projection. 
\end{definition}

\begin{definition} Let $i,j,k$ be integers, $1\leq i<j\leq n$, $i\leq
  k\leq n$. We let $\pi_{jik}\colon W\to \langle w_{jik}\rangle$ 
  be the canonical projection. 
\end{definition}

\begin{definition} Let $i$ be an integer, $1\leq i\leq n$. We let
  \[\Pi_i\colon V\to \langle v_{i1}, v_{i2},\ldots,
  v_{i(i-1)},v_{(i+1)i},\ldots,v_{ni}\rangle\]
be the canonical projection.
\end{definition}

The following observation follows immediately from the definition of
the $\varphi_i$:

\begin{lemma} For each $k$, $\varphi_k$ is one-to-one.
\end{lemma}

\begin{lemma} Let $\mathbf{w}\in\varphi_i(V)$. If
  $\pi_{rst}(\mathbf{w})\neq \mathbf{0}$, then $s\leq i\leq t$;
  moreover, at most one of the inequalities is strict.
\label{lemma:nonzerocoord}
\end{lemma}

\begin{lemma} Fix $n>1$, and let $i,j$ be integers satisfying $1\leq i<j\leq n$. Then
  $\varphi_i(V)\cap \varphi_j(V) = \{\mathbf{0}\}$.
\end{lemma}

\begin{proof} Assume $\varphi_i(\mathbf{v})\in \varphi_j(V)$, and
  that $\pi_{sr}(\mathbf{v})\neq \mathbf{0}$, where $1\leq r<s\leq
  n$. If $r\leq i$, then
  $\pi_{sri}(\varphi_i(\mathbf{v}))\neq\mathbf{0}$; but since
  $\varphi_i(\mathbf{v})\in\varphi_j(V)$, this implies that $r\leq
  j\leq i$, which is impossible.  If, on the other hand, $i<r$, then
  $\pi_{sir}(\varphi_i(\mathbf{v}))\neq \mathbf{0}$ and
  $\pi_{ris}(\varphi_i(\mathbf{v}))\neq\mathbf{0}$; again, since
  $\varphi_i(\mathbf{v})\in\varphi_j(V)$,
  Lemma~\ref{lemma:nonzerocoord} implies that $s=r=j$, which is also
  impossible.
\end{proof}

\begin{lemma} Fix $n>1$, $r\leq n$. Let $i_1,\ldots,i_r$ be integers,
  $1\leq i_1\leq \cdots\leq i_r\leq n$. Then
$\varphi_{i_1}(V) \cap \bigl\langle
\varphi_{i_2}(V),\ldots,\varphi_{i_r}(V)\bigr\rangle$
is of dimension $\binom{r-1}{2}$, with basis given by all vectors
$w_{ai_1b}-w_{bi_1a}$, $a,b\in \{i_2,\ldots,i_r\}$, $a> b$. A basis
for the pullback
\[ \varphi_{i_1}^{-1}\left(\Bigl\langle
\varphi_{i_2}(V),\ldots,\varphi_{i_r}(V)\Bigr\rangle\right)\]
is given by $\{v_{ab}\}$, with $a,b \in\{i_2,\ldots,i_r\}$, $a>b$.
\label{lemma:intersections}
\end{lemma}

\begin{proof} It is enough to   prove the last statement.
The result is trivial when $r=1,2$, so assume that
  $r\geq 3$. Certainly each $v_{ab}$ lies in the pullback, since
\[ \varphi_{i_1}(v_{ab}) = w_{ai_1b} -w_{bi_1a} =
\varphi_{b}(v_{ai_1}) - \varphi_a(v_{bi_i}).\]

If $\mathbf{w}\in\langle
\varphi_{i_2}(V),\ldots,\varphi_{i_r}(V)\rangle$, and
$\pi_{rst}(\mathbf{w})\neq \mathbf{0}$, then we must have $s\leq
i_j\leq t$ for some $j\in\{2,\ldots,r\}$, and with at most one
inequality strict. Let $\mathbf{v}\in V$ be a vector such that
$\varphi_i(\mathbf{v}) \in \langle
\varphi_{i_2}(V),\ldots,\varphi_{i_r}(V)\rangle$. Let $r,s$ be
integers, $1\leq r<s\leq n$, such that $\pi_{sr}(\mathbf{v})\neq
\mathbf{0}$. We want to show that $s,r\in \{i_2,\ldots,i_r\}$.

If $r\leq i_1$, then $\pi_{sri_1}(\varphi_i(\mathbf{v}))\neq
\mathbf{0}$, so we must have $r\leq i_j\leq i_1$ for some $j>1$, but
this is impossible, since $i_j>i_1$. Therefore, we must have
$r>i_1$. In that case, we have $\pi_{si_1r}(\mathbf{v})\neq\mathbf{0}$
and $\pi_{ri_1s}(\mathbf{v})\neq \mathbf{0}$. Thus we obtain that
there is some $j>1$ such that $i_1\leq i_j\leq r$ and at most one
inequality is strict, and there is some $k>1$ such that $i_1\leq i_k\leq
s$ and at most one inequality is strict. Since $i_1<i_j,i_k$, it
follows that $r=i_j$ and $s=i_k$, as desired.
\end{proof}

\begin{corollary} Fix $n>1$, $r\leq n$, and let $1\leq i_1\leq
  \cdots\leq i_r \leq n$ be
  integers. Then
\[ \dim\bigl(\langle
\varphi_{i_1}(V),\ldots,\varphi_{i_r}(V)\rangle\bigr) = r\binom{n}{2} -
\binom{r}{3}.\]
\end{corollary}

\begin{proof} Each $\varphi_k$ is injective, so we have:
\begin{eqnarray*}
\lefteqn{\dim\left(\langle
\varphi_{i_1}(V),\ldots,\varphi_{i_r}(V)\rangle\right)}\\
 &=&
\left(\sum_{k=1}^r \dim(\varphi_{i_k}(V))\right)
 -
\left(\sum_{k=1}^{r-2}\dim\Bigl(\varphi_{i_{r-k-1}}(V)\cap \langle
\varphi_{i_r}(V),\ldots,\varphi_{i_{r-k}}(V)\rangle \Bigr)\right)\\
& = & r\binom{n}{2} - \left(\sum_{k=1}^r \binom{k-1}{2}\right)
 =  r\binom{n}{2} - \binom{r}{3},
\end{eqnarray*}
as claimed.
\end{proof}

We have concentrated on the ``lowest index'' for simplicity. Of
course, given the definitions, our treatment has symmetry; for example:

\begin{prop} Let $i_1,\ldots,i_r\in\{1,\ldots,n\}$ be pairwise
  distinct. Then
\[ \varphi_{i_1}(V) \cap \Bigl\langle
\varphi_{i_2}(V),\ldots,\varphi_{i_r}(V)\Bigr\rangle\]
has dimension $\binom{r-1}{2}$. Moreover, a basis for the pullback is
given by the vectors $v_{ab}$, with $a,b\in\{i_2,\ldots,i_r\}$, and $a>b$.
\label{prop:countingwithsymmetry}
\end{prop}

\begin{proof} It is easy to verify that
\[ \varphi_{i_1}(V)\cap \Bigl\langle \varphi_{i_2}(V),\ldots,\varphi_{i_r}(V)\Bigr\rangle\]
is generated by the vectors $w_{abi_1}$ for
$a,b\in\{i_2,\ldots,i_r\}$, $a,i_1>b$, and the vectors
$w_{ai_1b}-w_{bi_1a}$ when $b>i_1$. These vectors are linearly
independent, and pulling them back gives the desired result.
\end{proof}

A somewhat different description of the intersections will be useful in the
following sections. 

\begin{definition} Fix $n>1$. We define $\Phi\colon V^n\to W$ to be
\[ \Phi(\mathbf{v}_1,\ldots,\mathbf{v}_n) = \varphi_1(\mathbf{v}_1) +
\cdots + \varphi_n(\mathbf{v}_n).\]
If there is danger of ambiguity, we use $\Phi_n$ to denote the map
corresponding to $n$.
\end{definition}

\begin{theorem} 
Fix $n>1$. The kernel of $\Phi$ is of dimension $\binom{n}{3}$. A basis for
${\rm ker}(\Phi)$ is as follows: for each choice of $a,b,c$, $1\leq
a<b<c\leq n$, the vector $(\mathbf{v}_1,\ldots,\mathbf{v}_n)\in V^n$
with
\[ \mathbf{v}_i = \left\{\begin{array}{ll}
v_{cb}&\mbox{if $i=a$;}\\
-v_{ca}&\mbox{if $i=b$;}\\
v_{ba}&\mbox{if $i=c$;}\\
\mathbf{0}&\mbox{otherwise.}
\end{array}\right.\]
\label{th:basisforkerPhi}
\end{theorem}

\begin{proof} Denote the element corresponding to $a,b,c$\/ by
  $\mathbf{v}_{(abc)}$. Note that each $\mathbf{v}_{(abc)}$ lies in
  ${\rm ker}(\Phi)$:
\[ \Phi(\mathbf{v}_{(abc)}) = \varphi_a(v_{cb}) + \varphi_b(-v_{ca}) +
  \varphi_c(v_{ba}) = (w_{cab} - w_{bac})  - w_{cab} + w_{bac} =
  \mathbf{0}.\]
Since $\Phi$ is surjective, $\dim(W) = n\dim(V) - \dim({\rm
  ker}(\Phi))$; therefore
\[ \dim\Bigl({\rm ker}(\Phi)\Bigr) = n\binom{n}{2} - 2\binom{n+1}{3} =
  \binom{n}{3}.\]
Thus, it is enough to show that the $\mathbf{v}_{(abc)}$ are linearly
independent. Assume that
$\sum\beta_{abc}\mathbf{v}_{(abc)}=\mathbf{0}$, where the sum is taken
over all triples of integers $a,b,c$ that satisfy $1\leq a<b<c\leq n$. Considering the $i$-th
component alone, we obtain
\[ \sum_{1\leq r<s<i}\!\!\!\! \beta_{rsi}v_{sr} -\!\!\!\! \sum_{1\leq r<i<s\leq
  n}\!\!\!\!\beta_{ris}v_{sr} +\!\!\!\! \sum_{i<r<s\leq n}\!\!\!\!\beta_{irs}v_{sr} =
  \mathbf{0}.\]
Thus, for each choice of $a,b,c$ with $i\in\{a,b,c\}$, we must have
  $\beta_{abc}=0$. This proves the given $\mathbf{v}_{(abc)}$ are
  linearly independent, and hence form a basis for $\ker(\Phi)$.
\end{proof}

\begin{notation} Fix $n>1$, and let $a,b,c$ be pairwise distinct
  integers, $1\leq a,b,c\leq n$. We will let $\mathbf{v}_{(abc)}$
  denote the element of $\ker(\Phi)$ described in the statement of
  Theorem~\ref{th:basisforkerPhi}.
\label{not:vabc}
\end{notation}

\begin{theorem}
Let $(\mathbf{v}_1,\ldots,\mathbf{v}_n)\in\ker(\Phi)$. Write
\[ \mathbf{v}_k = \sum_{1\leq i<j\leq n} \alpha_{ji}^{(k)}v_{ji},\]
\begin{itemize}
\item[(i)] If $i=k$ or $j=k$, then $\alpha_{ji}^{(k)}=0$; i.e.,
  $\Pi_k(\mathbf{v}_k) = \mathbf{0}$.
\item[(ii)] If $1\leq a<b<c\leq n$, then $\alpha_{ba}^{(c)} =
  \alpha_{cb}^{(a)} = -\alpha_{ca}^{(b)}$.
\item[(iii)] Fix $i,j$, $1\leq i<j\leq n$. Then
\begin{eqnarray*}
\Pi_{i}(\mathbf{v}_{j})& = &
\sum_{r=1}^{i-1}\left(-\alpha_{jr}^{(i)}\right)v_{ir}    + 
 \sum_{r=i+1}^{j-1} \alpha_{jr}^{(i)}v_{ri} +
\sum_{r=j+1}^n \left(-\alpha_{rj}^{(i)}\right)v_{ri},\\
\Pi_{j}(\mathbf{v}_i) & = &
\sum_{r=1}^{i-1}\left(-\alpha_{ir}^{(j)}\right)v_{jr} +
\sum_{r=i+1}^{j-1}\alpha_{ri}^{(j)}v_{jr} +
\sum_{r=j+1}^n\left(-\alpha_{ri}^{(j)}\right)v_{rj}.
\end{eqnarray*}
\end{itemize}
\label{th:descriptionkernelPhi}
\end{theorem}

\begin{proof} The first part follows either from
  Proposition~\ref{prop:countingwithsymmetry}, or from the description
  of the basis in Theorem~\ref{th:basisforkerPhi}. For part (ii), note
  that
\begin{eqnarray*}
\pi_{bac}\bigl(\varphi_1(\mathbf{v}_1) +\cdots + \varphi_n(\mathbf{v}_n)\bigr)
&=& \left(\alpha_{ba}^{(c)} - \alpha_{cb}^{(a)}\right)w_{bac},\\
\pi_{cab}\bigl(\varphi_1(\mathbf{v}_1)+\cdots +
\varphi_n(\mathbf{v}_n)\bigr) & = & \left(\alpha_{ca}^{(b)} +
\alpha_{cb}^{(a)}\right)w_{cab}.
\end{eqnarray*}
Since they must both be equal to zero, we obtain that
$\alpha_{ba}^{(c)}=\alpha_{cb}^{(a)}$ and $\alpha_{ca}^{(b)} =
-\alpha_{cb}^{(a)}$, as claimed.
Finally, for (iii), we know that
$\Pi_i(\mathbf{v}_i) = \Pi_j(\mathbf{v}_j)=\mathbf{0}$ from (i), so we can write:
\begin{eqnarray*}
\Pi_i(\mathbf{v}_j) & = & \sum_{r=1}^{i-1}\alpha_{ir}^{(j)}v_{ir} +
\sum_{r=i+1}^{j-1}\alpha_{ri}^{(j)}v_{ri} +
\sum_{r=j+1}^{n}\alpha_{ri}^{(j)}v_{ri},\\
\Pi_j(\mathbf{v}_i) & = & \sum_{r=1}^{i-1}\alpha_{jr}^{(i)}v_{jr} +
\sum_{r=i+1}^{j-1}\alpha_{jr}^{(i)}v_{jr} + \sum_{r=j+1}^n
\alpha_{rj}^{(i)}v_{rj},
\end{eqnarray*}
and applying (ii) gives the desired identities.
\end{proof}

\begin{corollary} Let $\mathbf{v}\in\ker(\Phi)$.
  $\Pi_j(\mathbf{v}_i)=\mathbf{0}$, then
  $\Pi_i(\mathbf{v}_j)=\mathbf{0}$. In particular, if
  $\mathbf{v}_i=\mathbf{0}$, then $\Pi_i(\mathbf{v}_j)=\mathbf{0}$ for
  all~$j$.
\end{corollary}

\begin{proof} If $\Pi_j(\mathbf{v}_i)=\mathbf{0}$, then
  $\alpha_{jr}^{(i)}=0$ for all $r$, so by
  Theorem~\ref{th:descriptionkernelPhi}(iii) it follows that
  $\Pi_i(\mathbf{v}_j)=\mathbf{0}$. The second assertion follows
  immediately.
\end{proof}

\begin{corollary} Let $\mathbf{v}\in\ker(\Phi)$, $\mathbf{v}\neq
  \mathbf{0}$. If 
$\mathbf{v} = (\mathbf{v}_1,\ldots,\mathbf{v}_n)$
then the dimension of $\langle
\mathbf{v}_1,\ldots,\mathbf{v}_n\rangle$ is at least $3$.
\label{cor:atleastthree}
\end{corollary}

\begin{proof} Write 
\[ \mathbf{v} = \!\!\!\!\!\sum_{1\leq a<b<c\leq n}\!\!\!\!\! \beta_{abc}\mathbf{v}_{(abc)}\]
Fix $a,b,c$ such that $1\leq a<b<c\leq n$, $\beta_{abc}\neq 0$. We
claim that $\mathbf{v}_a$, $\mathbf{v}_b$, and $\mathbf{v}_c$ are
linearly independent. Indeed, note that
$\Pi_a(\mathbf{v}_a)=\Pi_b(\mathbf{v}_b)=\Pi_c(\mathbf{v}_c) =
\mathbf{0}$, and $\pi_{cb}(\mathbf{v}_a) \neq\mathbf{0}$. Therefore,
if $\alpha_a\mathbf{v}_a + \alpha_b\mathbf{v}_b + \alpha_c\mathbf{v}_c
= \mathbf{0}$, then we must have $\alpha_a=0$. A symmetric argument
looking at $\pi_{ca}$ shows that $\alpha_b=0$, and considering
$\pi_{ba}$ shows that $\alpha_c=0$.
\end{proof}

In~\cite{capablep} we proved the following result by considering the
collection of images of a basis of~$X$ under the maps
$\varphi_1,\ldots,\varphi_n$, and showing they would necessarily be
linearly independent. We give a different proof here based on the
considerations above.

\begin{corollary}[Prop.~4.6 in~\cite{capablep}]
Fix $n>1$, and let $X$ be a subspace of~$V$. If $\dim(X)=1$, then
$\dim(X^*)=n$; if $\dim(X)=2$, then $\dim(X^*)=2n$.
\label{cor:maincountingarg}
\end{corollary}

\begin{proof} We prove the contrapositive.
Since $\dim(X^*) = n\dim(X) - \dim(X^n\cap\ker(\Phi))$, if
$\dim(X^*)<n\dim(X)$, then $X^n\cap\ker(\Phi)\neq\{\mathbf{0}\}$. 

Let $\mathbf{v}=(\mathbf{v}_1,\ldots,\mathbf{v}_n)\in
X^n\cap\ker(\Phi)$, $\mathbf{v}\neq \mathbf{0}$. Then $\mathbf{v}_i\in
X$ for $i=1,\ldots,n$, so by Corollary~\ref{cor:atleastthree},
$\dim(X)\geq 3$, as claimed.
\end{proof}

As mill become apparent in Section~\ref{sec:nequalfour}, it would be
useful to make Corollary~\ref{cor:atleastthree} somewhat more
precise. Unfortunately it appears we cannot do this easily, as the
next sequence of results shows. We will show that for any $k$, $3\leq
k\leq n$, $k\neq 4$, there exists
$(\mathbf{v}_1,\ldots,\mathbf{v}_n)\in\ker(\Phi_n)$ such that the span
of $\mathbf{v}_1,\ldots,\mathbf{v}_n$ is $k$-dimensional. I do not
know whether one can also obtain a $4$-dimensional span. 

\begin{lemma} Fix $n>1$, and let $k$ be an integer, $3\leq k\leq
  n$. If there exists
  $(\mathbf{v}_1,\ldots,\mathbf{v}_n)\in\ker(\Phi_n)$ such that
  $\dim\bigl(\langle \mathbf{v}_1,\ldots,\mathbf{v}_n\rangle\bigr) =
  k$, then for any $m\geq n$ there exists
  $(\mathbf{w}_1,\ldots,\mathbf{w}_m)\in\ker(\Phi_m)$ such that
  $\dim\bigl(\langle\mathbf{w}_1,\ldots,\mathbf{w}_m\rangle\bigr)=k$.
\label{lemma:increasen}
\end{lemma}

\begin{proof}
We can embedd $\ker(\Phi_n)$ into $\ker(\Phi_{n+1})$
by mapping $(\mathbf{v}_1,\ldots,\mathbf{v}_n)\in\ker(\Phi_n)$
to $(\mathbf{v}_1,\ldots,\mathbf{v}_n,\mathbf{0})$. The dimension of
the span of the components is clearly unaffected by this embedding.
\end{proof}

\begin{lemma} Fix $n>1$. If
  $(\mathbf{v}_1,\ldots,\mathbf{v}_n)\in\ker(\Phi_n)$, then 
\[\mathbf{w}=(\mathbf{v}_1,\ldots,\mathbf{v}_n,\mathbf{0},\mathbf{0},\mathbf{0})+\mathbf{v}_{(abc)}\in\ker(\Phi_{n+3}),\]
where $a=n+1$, $b=n+2$, and $c=n+3$ (using the notation described
in~\ref{not:vabc}), and if we let
$\mathbf{w}=(\mathbf{w}_1,\ldots,\mathbf{w}_{n+3}$, then
\[
\dim\bigl(\langle\mathbf{w}_1,\ldots,\mathbf{w}_{n+3}\rangle\bigr)
= \dim\bigl(\langle\mathbf{v}_1,\ldots,\mathbf{v}_n\rangle\bigr)+3.\]
\label{lemma:addthreetodim}
\end{lemma}

\begin{proof} The first part follows by embedding $\ker(\Phi_n)$ into
  $\ker(\Phi_{n+3})$ by appending three zeros; then both vectors lie
  in $\ker(\Phi_{n+3})$, hence so does their sum $\mathbf{w})$. For
  the second part, simply note that $\mathbf{w}_i=\mathbf{v}_i$ for
  $1\leq i\leq n$, and $\mathbf{w}_{n+1}=v_{n+3,n+2}$,
  $\mathbf{w}_{n+2}=-v_{n+3,n+1}$, and
  $\mathbf{w}_{n+3}=v_{n+2,n+1}$. This adds three to the dimension of
  the span of $\mathbf{v}_1,\ldots,\mathbf{v}_n$, since
  $\Pi_k(\mathbf{v}_i)=\mathbf{0}$ for $1\leq i\leq n$, $k>n$.
\end{proof}

\begin{lemma} If $k=3\ell$ or $k=3\ell+2$, $\ell\geq 1$, and $n\geq k$, then there
  exists $(\mathbf{v}_1,\ldots,\mathbf{v}_n)\in\ker(\Phi_n)$ such that
$\dim\bigl(\langle\mathbf{v}_1,\ldots,\mathbf{v}_n\rangle\bigr) =
k$.
The result also holds if $k=3\ell+1$ and $\ell\geq 2$.
\label{lemma:mostvaluesofk}
\end{lemma}
\begin{proof} For $k=3$, any nontrivial element of $\ker(\Phi_3)$ will do.
  Then we can apply Lemmas~\ref{lemma:increasen}
  and~\ref{lemma:addthreetodim} to obtain the result whenever $k$ is a
  multiple of~$3$.
For $k=5$ and $n=5$, the element 
$\mathbf{v}_{(123)} +
  \mathbf{v}_{(145)} = \bigl( v_{32}+v_{54}, -v_{31}, v_{21}, -v_{51},
  v_{41}\bigr) \in \ker(\Phi_5)$
gives an example of dimension~$5$, and again applying Lemmas~\ref{lemma:increasen}
  and~\ref{lemma:addthreetodim} we obtain the result whenever $k\equiv
  2\pmod{3}$, $k\geq 5$. 

If $k=7$, then we have
$\mathbf{v}_{(123)} + \mathbf{v}_{(145)} +
\mathbf{v}_{(176)}\in\ker(\Phi_7)$
which gives an example of dimension~$7$, from which we can obtain any $k\equiv
1\pmod{3}$, $k\geq 7$ as above.
\end{proof}

I currently do not know if one can find an element of $\ker(\Phi_n)$
which will yield a subspace of dimension exactly four for some~$n$.
We can, however, show that it is impossible if  $n=4$:

\begin{prop} If
  $(\mathbf{v}_1,\mathbf{v}_2,\mathbf{v}_3,\mathbf{v}_4)$ is
  nontrivial and lies in $\ker(\Phi_4)$, then 
\[\dim(\langle
  \mathbf{v}_1,\mathbf{v}_2,\mathbf{v}_3,\mathbf{v}_4\rangle)=3.\]
\label{prop:kerPhifornfour}
\end{prop}

\begin{proof}
From Theorem~\ref{th:descriptionkernelPhi} it follows that we can
write the vectors $\mathbf{v}_i$ as follows:
\begin{eqnarray*}
\mathbf{v}_1 & = & \beta_{123} v_{32} + \beta_{124} v_{42} + \beta_{134} v_{43},\\
\mathbf{v}_2 & = & -\beta_{123}v_{31} - \beta_{124} v_{41} + \beta_{234} v_{43},\\
\mathbf{v}_3 & = & \beta_{123} v_{21} - \beta_{134} v_{41} - \beta_{234} v_{42},\\
\mathbf{v}_4 & = & \beta_{124} v_{21} + \beta_{134} v_{31} + \beta_{234} v_{32},
\end{eqnarray*}
for some choice of coefficients $\beta_{123}$, $\beta_{124}$,
$\beta_{134}$, and~$\beta_{234}$. If not all $\beta_{abc}$ are equal
to zero, then we have the following nontrivial linear relation between
the four vectors:
\[ \beta_{234}\mathbf{v}_1 +\beta_{124}\mathbf{v}_3 =  \beta_{134}\mathbf{v}_2 +
 \beta_{123}\mathbf{v}_4.\]
By Corollary~\ref{cor:atleastthree} the four
vectors span a space of dimension at least three. Thus, the subspace
is of dimension exactly three, and if $\beta_{abc}\neq 0$, then a basis for the
subspace is given by $\mathbf{v}_a$, $\mathbf{v}_b$,
and~$\mathbf{v}_c$. This establishes the result for $n=4$.
\end{proof}

\section{Dimension counting.}\label{sec:dimcounting}

In this section we will obtain bounds for $\dim(X^*)$ in terms of
$\dim(X)$. To see why this is interesting, consider the following two observations:

\begin{prop}
Let $X<V$. Assume that for all subspaces $Y$ of\/ $V$, if\/ $Y$ properly contains $X$
then $Y^*$ properly contains $X^*$. Then $X=X^{**}$.
\label{prop:strictlylarger}
\end{prop}

\begin{proof} If $X^{**}$ properly contains $X$, then $X^{***}$ would
  properly contain $X^*$. But $X^{***}=X^*$, a contradiction.
\end{proof}

\begin{corollary}
Fix $n>1$, and suppose that $f(k)$ is a function such that for all
$k$-dimensional subspaces of~$V$, $\dim(X^*)\geq nk-f(k)$. 
If $f(k+1)<n$, then all
subspaces of dimension $k$ are closed. 
\end{corollary}

\begin{proof} It is trivial that $\dim(X^*)\leq n\dim(X)$.
If $f(k+1)<n$ then for subspaces $X$ and~$Y$ with
  $\dim(X)=k$ and $\dim(Y)=k+1$ we have
\[ \dim(Y^*) \geq n(k+1) - f(k+1) > nk \geq \dim(X^*),\]
and thus, by Proposition~\ref{prop:strictlylarger}, it follows that
all $X$ of dimension $k$ are closed.
\end{proof}

The smallest possible value for $f(k)$ is given by
\[ f(k) = \max\bigl\{ \dim(X^n\cap\ker(\Phi))\,\bigm|\,
\dim(X)=k\}.\]
Our objective in this section is to find an expression
for $f(k)$ in terms of~$k$.  The main workhorse in our calculations
will Lemma~\ref{lemma:basicoverlapbound} below.  The idea is to find
the value of $\dim(X^n\cap\ker(\Phi))$ by examining the ``partial
intersections''; namely, the intersections of the form
\[ \ker(\Phi)\cap \Bigl\langle
(\mathbf{0},\ldots\mathbf{0},\mathbf{v}_i,
\mathbf{v}_{i+1},\ldots,\mathbf{v}_n) \,\Bigm|\, \mathbf{v}_j\in
X\Bigr\rangle,\]
as $i$ ranges from $1$ to $n-2$ (when $i=n-1$ or $i=n$, the
intersection is trivial by Corollary~\ref{cor:atleastthree}). For a fixed $i$, we can consider the
subspace of $X$ consisting of all vectors $\mathbf{v}_i$ which can be
``completed'' to an element of $\ker(\Phi)$ by appending $\mathbf{0}$
prior to it, and any vector of~$X$ after. This is the same as
considering the pullbacks
$ X \cap
\varphi_i^{-1}\left(\langle\varphi_{i+1}(X),\ldots,\varphi_{n}(X)\rangle\right)$.
It is easy to verify that the sum of the dimensions of these pullbacks
is equal to the dimension of $X^n\cap\ker(\Phi)$.  We will first use
the dimension of these pullbacks to establish a lower bound for the
dimension of~$X$; then we will turn around and use these calculations
to give an upper bound for the dimension of the pullbacks in terms of
the dimension of~$X$.

Making the bounds as precise as possible, however, requires one to
keep track of a lot of information; this in turn requires the use of
multiple indices and subindices in the proof, for which I apologize in
advance. To illustrate the ideas and help the reader navigate through
the proof, we will first present an example. This is not an example in
the sense of a specific $X$, but rather a run-through the main part of
the analysis we will perform with specific values for the indices and
some of the variables.

\begin{example} Set $n=6$, and let $X$ be a subspace of~$V$. We will
be interested in bounding above the dimension of $Z_i$ in terms of
$\dim(X)$, where
\[Z_i = X \cap \varphi_i^{-1}\Bigl(\langle
\varphi_{i+1}(X),\ldots,\varphi_6(X)\rangle\Bigr);\]
i.e., $Z_i$ consists of all $\mathbf{v}\in X$ for which there exist
$\mathbf{v}_{i+1},\ldots,\mathbf{v}_6$ in~$X$ such that
\[
(\mathbf{0},\ldots,\mathbf{0},\mathbf{v},\mathbf{v}_{i+1},\ldots,\mathbf{v}_6)
  \in X^6\cap\ker(\Phi).\]
To do this, we will obtain a lower bound for $\dim(X)$ in terms of $\dim(Z_i)$.
To further fix ideas, set $i=2$.
Note that by Lemma~\ref{lemma:intersections} (or
  Theorem~\ref{th:descriptionkernelPhi}) we must have
  $\Pi_1(Z_2)=\Pi_2(Z_2)=\mathbf{0}$.
Assume that the dimension of $Z_2$ is $4$. Order all pairs $(j,i)$
  lexicographically from right to left, so $(j,i)<(b,a)$ if and only
  if $i<a$, or $i=a$ and $j<b$. Then considering all pairs in order,
  find out which pairs $(b,a)$ have $\pi_{ba}(Z_2)\neq \mathbf{0}$. In
  this example, say that it is all possible pairs:
\[ (4,3), (5,3), (6,3), (5,4), (6,4), (6,5).\]
Doing row reduction, we can find a basis $\mathbf{v}_{1,2}$,
$\mathbf{v}_{2,2}$, $\mathbf{v}_{3,2}$, and~$\mathbf{v}_{4,2}$
for $Z_2$ (the second index refers to the fact that
these vectors are in the second component of an element of
$\ker(\Phi)$), satisfying that the ``leading pair'' (smallest nonzero
component) of each is strictly smaller than that of its successors, 
and all other vectors have zero component for that pair. For example, 
\[ \begin{array}{rclcrcl}
\mathbf{v}_{1,2} & = & v_{43} + \alpha_1 v_{53} + \alpha_2 v_{64},&\qquad&
\mathbf{v}_{3,2} & = & v_{54} + \gamma v_{64},\\
\mathbf{v}_{2,2} & = & v_{63} + \beta v_{64},&\qquad&
\mathbf{v}_{4,2} & = & v_{65},
\end{array}\]
for some coefficients $\alpha_1,\alpha_2,\beta,\gamma\in\mathbb{F}_p$.
We know there exist vectors $\mathbf{v}_{i,3}$, $\mathbf{v}_{i,4}$,
$\mathbf{v}_{i,5}$, $\mathbf{v}_{i,6}$ such that
$
(\mathbf{0},\mathbf{v}_{i,2},\mathbf{v}_{i,3},\mathbf{v}_{i,4},\mathbf{v}_{i,5},\mathbf{v}_{i,6})\in
X^6\cap\ker(\Phi)$ for $i=1,2,3,4$. Naturally, $X$ contains all
twenty vectors, but there will normally be some linear
dependencies between them: some may even be equal to~$\mathbf{0}$. We
want to extract a subset that we know is linearly independent in some
systematic fashion. First let us consider the information we can obtain about
these vectors from our knowledge of the vectors~$\mathbf{v}_{i,2}$.

Since
$(\mathbf{0},\mathbf{v}_{i,2},\mathbf{v}_{i,3},\mathbf{v}_{i,4},\mathbf{v}_{i,5},\mathbf{v}_{i,6})$
lies in $\ker(\Phi)$, we can use
Theorem~\ref{th:descriptionkernelPhi}(iii) to describe the $\Pi_i$-image of each
vector $\mathbf{v}_{i,j}$, where $i\leq 2$ and $j>2$.
The $\Pi_1$-image must be
trivial, and for the $\Pi_2$ image we obtain
the following:
\begin{equation*}
\begin{aligned}
\Pi_2(\mathbf{v}_{1,3}) & = v_{42} + \alpha_1 v_{52},\\
\Pi_2(\mathbf{v}_{1,4}) & = -v_{32} + \alpha_2 v_{62},\\
\Pi_2(\mathbf{v}_{1,5}) & = -\alpha_1 v_{32},\\
\Pi_2(\mathbf{v}_{1,6}) & = -\alpha_2 v_{42}.
\end{aligned}
\qquad\qquad\qquad
\begin{aligned}
\Pi_2(\mathbf{v}_{2,3}) & = v_{62},\\
\Pi_2(\mathbf{v}_{2,4}) & = \beta v_{62},\\
\Pi_2(\mathbf{v}_{2,5}) & = \mathbf{0},\\
\Pi_2(\mathbf{v}_{2,6}) & = -v_{32} - \beta v_{42}.
\end{aligned}
\end{equation*}
\begin{equation*}
\begin{aligned}
\Pi_2(\mathbf{v}_{3,3}) & = \mathbf{0},\\
\Pi_2(\mathbf{v}_{3,4}) & = v_{52} + \gamma v_{62},\\
\Pi_2(\mathbf{v}_{3,5}) & = - v_{42},\\
\Pi_2(\mathbf{v}_{3,6}) & = -\gamma v_{42}.
\end{aligned}
\qquad\qquad\qquad
\begin{aligned}
\Pi_2(\mathbf{v}_{4,3}) & = \mathbf{0},\\
\Pi_2(\mathbf{v}_{4,4}) & = \mathbf{0},\\
\Pi_2(\mathbf{v}_{4,5}) & = v_{62},\\
\Pi_2(\mathbf{v}_{4,6}) & = -v_{52}.
\end{aligned}
\end{equation*}

One way to obtain these without too much confusion is as follows: to
find $\Pi_2\left(\mathbf{v}_{j,k}\right)$, go through the expression
for $\mathbf{v}_{j,2}$ replacing all indices $k$ by $2$, remembering that
$v_{ab}=-v_{ba}$. Any $v_{ba}$ in which neither $a$ nor~$b$ are equal
to $k$ are simply removed. 

To systematically extract from this list a set of linearly
independent vectors, we proceed as follows: consider all the pairs
which are leading components of the basis vectors; in this case,
$(4,3)$, $(6,3)$, $(5,4)$, and~$(6,5)$.
The individual
indices that occur are $3$, $4$, $5$, and $6$. For each of them, we
identify the smallest pair in which it occurs. Thus, $3$ first occurs
in pair number one, as does $4$. The index $5$ first occurs in pair
number three, and $6$ first occurs in pair number two.

Since the first pair in which $3$ appears is the \textit{first} pair,
$(4,3)$, if we consider the vector we obtain when we replace the index
$4$ (the other index in the pair we found for $3$) from the vector
which has $(4,3)$ as its leading component, we will obtain a vector whose
first nontrivial component is $(3,2)$. That is, the vector
$\mathbf{v}_{1,4}$ (replacing $4$ in the \textit{first} vector).
The next index is~$4$, again in the \textit{first} pair. If we replace
the other index, $3$, in the \textit{first} vector, i.e. when we look
at $\mathbf{v}_{1,3}$, we obtain a vector
with nontrivial $(4,2)$, component, and for which all
$(j,i)$-components with $(j,i)<(4,2)$ are trivial. 
Then take the index~$5$: it first occurs in the \textit{third} pair,
paired with~$4$, so the vector we obtain by replacing $4$ in the
\textit{third} vector, i.e.~the vector $\mathbf{v}_{3,4}$, is a
vector with nontrivial $(5,2)$ component and trivial $(j,i)$
component for all $(j,i)<(5,2)$. For the index $6$ we will take
$\mathbf{v}_{2,3}$ (since $6$ first occurs in the second pair,  paired with~$3$)
which gives a vector with nontrivial $(6,2)$ component and trivial
$(j,i)$ component for all $(j,i)<(6,2)$. 
In summary, we want to consider, in addition to the basis
for $Z_2$, the vectors
$\mathbf{v}_{1,4}, \mathbf{v}_{1,3}, \mathbf{v}_{3,4},
\mathbf{v}_{2,3}$
corresponding, respectively, to the indices $3$, $4$, $5$,
and~$6$. The choices we have made ensure that the $\Pi_2$-images of
these vectors are linearly independent, 
and so the vectors themselves must be linearly independent.
Since $\Pi_2(Z_2)=\mathbf{0}$, the subset formed of the basis for $Z_2$
together with these four vectors is a linearly independent subset
of~$X$; so we can conclude that~$X$ must have dimension at least
$8$. What is more, note that none of these last four vectors will
occur in a similar analysis involving $Z_3$: when performing a similar
analysis, all vectors will have trivial $\Pi_i$-image when $i<3$. 
Note as well that the number of indices, in this case $4$, must
satisfy $\dim(Z_2)\leq\binom{4}{2}$, since we need to be able to
obtain at least $\dim(Z_2)$ pairs out of the indices that occur.~$\Box$
\end{example}

What ensures that this process will work the way we want is how
we choose the vectors of the basis and the vectors that ``correspond''
to each index. The former count towards the value of
$\dim(X^n\cap\ker(\Phi))$, while the latter may be removed from
consideration when we move on to $Z_{i+1}$. 
This is all done in generality in the proof
of the following promised lemma:

\begin{lemma} Let $X$ be a subspace of $V$. For each $i$, $1\leq i\leq
  n$, let 
\[Z_i = X \cap \varphi_i^{-1}\Bigl(\langle
\varphi_{i+1}(X),\ldots,\varphi_n(X)\rangle\Bigr);\]
i.e., $Z_i$ consists of all $\mathbf{v}\in X$ for which there exist
$\mathbf{v}_{i+1},\ldots,\mathbf{v}_n$ in~$X$ such that
\[
(\mathbf{0},\ldots,\mathbf{0},\mathbf{v},\mathbf{v}_{i+1},\ldots,\mathbf{v}_n)
  \in X^n\cap\ker(\Phi).\]
If
$\dim\bigl(X\cap\langle v_{sr}\,|\,i\leq r<s\leq
  n\rangle\bigr)=d_i$ and
$\dim(Z_i)=r_i$, then $r_i\leq \binom{d_i-r_i}{2}$. Morevoer, if
$s_i$ is the smallest positive integer such that
$r_i\leq\binom{s_i}{2}$, then we must have $d_{i+1}\leq d_i-s_i$.
\label{lemma:basicoverlapbound}
\end{lemma}

\begin{proof} Fix $i_0$, $1\leq i_0\leq n$. For simplicity, write $r=r_{i_0}$.
By Theorem~\ref{th:descriptionkernelPhi}, if $\mathbf{v}\in Z_{i_0}$ then
$\Pi_i(\mathbf{v})=\mathbf{0}$ for all $i\leq i_0$. 

Let $\mathbf{v}_{1i_0},\ldots,\mathbf{v}_{ri_0}$ be a basis for
 $Z_{i_0}$. We will modify it as follows:

Order all pairs $(j,i)$, $i_0<i<j\leq n$ by letting $(j,i)<(b,a)$ if
 and only if $i<a$ or $i=a$ and $j<b$ (lexicographically from right to
 left).  Let $(j_1,i_1)$ be the smallest pair for which
 $\pi_{j_1i_1}(\mathbf{v}_{ki_0})\neq\mathbf{0}$ for some $k$, $1\leq
 k\leq r$. Reordering if necessary we may assume $k=1$. Replacing
 $\mathbf{v}_{1i_0}$ with a scalar multiple of itself and adding
 adequate multiples to the remaining $\mathbf{v}_{ki_0}$ if necessary
 we may also assume that
\[ \pi_{j_1i_1}\left(\mathbf{v}_{ki_0}\right) = \left\{\begin{array}{ll}
v_{j_1i_1}&\mbox{if $k=1$;}\\
\mathbf{0}&\mbox{if $k\neq 1$.}
\end{array}\right.\]

Let $(j_2,i_2)$ be the smallest pair for which
$\pi_{j_2i_2}(\mathbf{v}_{ki_0})\neq\mathbf{0}$ for some $k$, $2\leq
k\leq r$. Again we may assume $k=2$, and that 
\[\pi_{j_2i_2}\left(\mathbf{v}_{ki_0}\right) = \left\{\begin{array}{ll}
v_{j_2i_2}&\mbox{if $k=2$;}\\
\mathbf{0}&\mbox{if $k\neq 2$.}
\end{array}\right.\]

Proceeding in the same way for $k=3,\ldots,r$, we obtain an ordered
list of pairs
$(j_1,i_1)<(j_2,i_2)<\ldots<(j_r,i_r)$ and a basis
$\mathbf{v}_{1i_0},\ldots,\mathbf{v}_{ri_0}$ such that
\[ \pi_{j_{\ell}i_{\ell}}\left(\mathbf{v}_{ki_0}\right) =
  \left\{\begin{array}{ll}
v_{j_{\ell}i_{\ell}}&\mbox{if $\ell=k$,}\\
\mathbf{0}&\mbox{if $\ell\neq k$;}
\end{array}\right.\]
and such that
$\pi_{ba}\left(\mathbf{v}_{ki_0}\right) = \mathbf{0}\quad\mbox{for
  all $(b,a)<(j_k,i_k)$}$.
Write $\displaystyle \mathbf{v}_{ki_0} =\!\!\!\!\! \sum_{i_0<i<j\leq n}
\!\!\!\!\!\alpha_{ji}^{(k,i_0)}v_{ji}$.  From the above we have:
\[ \alpha_{ji}^{(k,i_0)} = \left\{ \begin{array}{ll}
1 &\mbox{if $(j,i)=(j_k,i_k)$,}\\
0 &\mbox{if $(j,i)<(j_k,i_k)$.}
\end{array}\right.\]
For $k=1,\ldots,r$ and $i=i_0+1,\ldots,n$, let $\mathbf{v}_{ki}$ be
vectors in~$X$ such that
\[
\left(\mathbf{0},\ldots,\mathbf{0},\mathbf{v}_{ki_0},\mathbf{v}_{ki_0+1},\ldots,\mathbf{v}_{kn}\right)\in
\ker(\Phi)\cap X^n.\]
By Theorem~\ref{th:descriptionkernelPhi}(iii) we have
\begin{equation*}
\Pi_{i_0}\left(\mathbf{v}_{kj}\right) = \sum_{m=i_0+1}^{j-1}
\!\!\!\alpha_{jm}^{(k,i_0)}v_{mi_0} - \sum_{m=j+1}^n
\!\!\!\alpha_{mj}^{(k,i_0)}v_{mi_0}.
\end{equation*}
For simplicity, set $\alpha_{ji}^{(k,i_0)} = -\alpha_{ij}^{(k,i_0)}$,
and $\alpha_{jj}^{(k,i_0)}=0$;
then we can rewrite the above expression as:
\begin{equation}
\Pi_{i_0}\left(\mathbf{v}_{kj}\right) = \sum_{m=i_0+1}^{n}
\!\!\!\alpha_{jm}^{(k,i_0)}v_{mi_0}.
\label{eq:generalformula}
\end{equation}

Let $s$ be the cardinality of the set $\{i_1,j_1,\ldots,i_r,j_r\}$;
that is, $s$ is the number of distinct indices that occur in the list
$(j_1,i_1),\ldots,(j_r,i_r)$. Note that $r\leq \binom{s}{2}$. Let
$a_1<a_2<\cdots < a_s$ be the list of these distinct indices. For each
$\ell$ with $1\leq \ell\leq s$, let $(j_{k(\ell)},i_{k(\ell)})$ be the
smallest pair among $(j_1,i_1),\ldots,(j_r,i_r)$ that has
$a_{\ell}\in\{ i_{k(\ell)}, j_{k(\ell)}\}$. If 
$a_{\ell}=i_{k(\ell)}$, let $b_{\ell} = j_{k(\ell)}$; if
$a_{\ell}=j_{k(\ell)}$, let $b_{\ell}=i_{k(\ell)}$. Consider the
following list of vectors from~$X$:
\[\mathbf{v}_{1i_0},\mathbf{v}_{2i_0},\ldots,
\mathbf{v}_{ri_0},\mathbf{v}_{k(1)b_1},
\mathbf{v}_{k(2)b_2},\ldots,\mathbf{v}_{k(s)b_s}.\] 
Note that all of these vectors lie in $X\cap\langle v_{ji} \,|\,
i_0\leq i<j\leq n\rangle$.
We will show that these vectors are linearly independent. 
Since $\mathbf{v}_{1i_0},\ldots,\mathbf{v}_{ri_0}$ are
linearly independent and $\Pi_{i_0}(\mathbf{v}_{ki_0})=\mathbf{0}$ for
$k=1,\ldots,r$, it suffices to show that
$\Pi_{i_0}(\mathbf{v}_{k(1)b_1}),\ldots,\Pi_{i_0}(\mathbf{v}_{k(s)b_s})$
are linearly independent.

First, from $(\ref{eq:generalformula})$ we have
$ \pi_{a_{\ell}i_0}\left(\mathbf{v}_{k(m)b_m}\right) =
\alpha_{b_ma_{\ell}}^{(k(m),i_0)}$. We claim that if $\ell<m$, then
$\alpha_{b_m a_{\ell}}^{(k(m),i_0)}=0$. By construction, this claim
will follow if we can show that either $a_{\ell}=b_m$, or else the
pair made up of $b_m$ and $a_{\ell}$ is strictly smaller than the pair
made up of $a_m$ and $b_m$ (which is equal to $(j_{k(m)},i_{k(m)})$);
the claim will then follow because $\alpha_{ba}^{(k,i_0)}=0$ whenever
$(b,a)<(j_k,i_k)$.  Indeed, we know that $a_{\ell}<a_m$. If
$a_m=i_{k(m)}$ and $b_m=j_{k(m)}$, then replacing $a_m$ in the pair
$(b_m,a_m)$ with something smaller (namely $a_{\ell}$) gives a smaller
pair: $(b_m,a_{\ell})<(b_m,a_m)$. If, on the other hand, we have
$a_m=j_{k(m)}$ and $b_m=i_{k(m)}$, then if $a_{\ell}>b_m$ we have
$(a_{\ell},b_m)<(a_m,b_m)$, and if $a_{\ell}<b_m$ then we also have
$(b_m,a_{\ell})<(a_m,b_m)$. The only remaining possibility is
$a_{\ell}=b_m$, which is of course no trouble.

Thus, we conclude that
$\alpha_{b_ma_{\ell}}^{k(m),i_0}=0$ whenever $\ell<m$. To 
see that the vectors
$\Pi_2(\mathbf{v}_{k(1)b_1}),\ldots,\Pi_2(\mathbf{v}_{k(s)b_s})$ are
linearly independent, note that
\[ \pi_{a_{\ell}i_0}\left(\mathbf{v}_{k(m)b_m}\right) =
\alpha_{b_ma_{\ell}}^{(k(m),i_0)} = \left\{\begin{array}{ll}
\mathbf{0}&\mbox{if $\ell<m$,}\\
v_{b_{\ell}a_{\ell}}&\mbox{if $m=\ell$.}
\end{array}\right.\]
Therefore, if $\beta_1\Pi_{i_0}(\mathbf{v}_{k(1)b_1}) + \cdots +
\beta_s\Pi_{i_0}(\mathbf{v}_{k(s)b_s})=\mathbf{0}$, then $\beta_1=0$ since
the only vector with nontrivial $(a_1,i_0)$-component is
$\Pi_{i_0}(\mathbf{v}_{k(1)b_1})$. Hence $\beta_2=0$, because the only
remaining vector with nontrivial $(a_2,i_0)$-component is
$\Pi_{i_0}(\mathbf{v}_{k(2)b_2})$; and continuing this way we conclude
$\beta_j=0$ for all $j$. So
the vectors are indeed linearly independent. Thus we have established
that
\[
\mathbf{v}_{1i_0},\mathbf{v}_{2i_0},\ldots,\mathbf{v}_{ri_0},\mathbf{v}_{k(1)b_1},\mathbf{v}_{k(2)b_2},\ldots,\mathbf{v}_{k(s)b_s}\]
is a collection of linearly independent vectors in~$X\cap\langle
v_{sr}\,|\,i_0\leq r<s\leq n\rangle$. 

Thus we conclude that $d_{i_0}\geq r+s$. Since $r\leq
\binom{s}{2}$, it follows that 
\[r\leq \binom{s}{2} \leq \binom{d_{i_0}-r}{2},\]
as claimed.

To complete the proof, it only remains to establish the upper bound on
$d_{i_0+1}$. We have
$d_{i_0} = d_{i_0+1} + \dim\bigl( \langle v_{ji}\,|\, i_0\leq
i<j\leq n\rangle\cap\{\mathbf{v}\in X\,|\, \Pi_{i_0}(\mathbf{v})\neq \mathbf{0}\}\bigr)$.
Since the vectors $\mathbf{v}_{k(1)b_1},\ldots,\mathbf{v}_{k(s)b_s}$ are linearly
independent, have nontrivial $\Pi_{i_0}$ projection, and lie in
$X \cap \bigl\langle
v_{ji}\,\bigm|\,i_0\leq i<j\leq n\bigr\rangle$,
we have $d_{i_0} \geq d_{i_0+1}+s$. Moreover, since $r\leq
\binom{s}{2}$, we also have $s_{i_0}\leq s$; therefore,
$d_{i_0+1} \leq d_{i_0}-s \leq d_{i_0}-s_{i_0}$,
as desired.
\end{proof}

Note that $Z_{n-1}$ and $Z_n$ are always trivial. 

\begin{definition} Let $d$ be a nonnegative integer. We define $r(d)$
  to be the largest integer such that $r(d)\leq d$ and $r(d)\leq
  \binom{d-r(d)}{2}$.
\end{definition}

\begin{theorem}
Fix $n>1$ and let $X<V$. Fix $i_0$, $1\leq i_0\leq n-2$, and let 
\[ Z_{i_0} = X\cap \varphi_{i_0}^{-1}\left(\Bigl\langle
\varphi_{i_0+1}(X),\ldots,\varphi_n(X)\Bigr\rangle\right).\]
If $\dim(X\cap \langle v_{ji}\,|\, i_0\leq i<j\leq n\rangle) = d$,
then $\dim(Z_{i_0})\leq r(d)$.
Equivalently,
\begin{equation}
\dim(Z_{i_0})\leq d -
\left\lceil\frac{\sqrt{8d+1}-1}{2}\right\rceil
\label{eq:maxoverlap}
\end{equation}
where $\lceil x\rceil$ is the smallest integer greater than or equal
to $x$.
\label{th:overlaplevelktermsoft}
\end{theorem}

\begin{proof} Let $\dim(Z_{i_0})=r$. By
  Lemma~\ref{lemma:basicoverlapbound}, $r\leq\binom{d-r}{2}$, so
  $r\leq r(d)$, as claimed. From $r(d)\leq \binom{d-r(d)}{2}$ we
  easily obtain $(\ref{eq:maxoverlap})$.
\end{proof}

We have two other ways of describing the function $r(d)$, which will
prove useful below:

\begin{corollary}
Let $d$ be a positive integer. Then $r(d)$ is the number of
nontriangular numbers strictly less than $d$.
 Equivalently, if we write
$d = \binom{t}{2} + s$, with $0< s \leq t$,
then $r(d) = \binom{t-1}{2} + (s-1)$.
\label{cor:overlapleveltriangular}
\end{corollary}

\begin{proof} Since $r(d)\leq \binom{d-r(d)}{2} \leq
  \binom{(d+1)-r(d)}{2}$, it follows that $r(d+1)\geq r(d)$. 
We also have
\[ r(d)+2 > r(d)+1 > \binom{d-(r(d)+1)}{2} =
\binom{(d+1)-(r(d)+2)}{2},\]
so $r(d+1)<r(d)+2$.
If $r(d)<\binom{d-r(d)}{2}$, then
\[ r(d)+1\leq \binom{d-r(d)}{2} = \binom{(d+1)-(r(d)+1)}{2},\]
so $r(d+1)\geq r(d)+1$ and in this case we have $r(d+1)=r(d)+1$. If
$r(d)=\binom{d-r(d)}{2}$, then $r(d)+1>\binom{(d+1)-(r(d)+1)}{2}$,
hence $r(d+1)<r(d)+1$ and we conclude that $r(d+1)=r(d)$.  In summary,
we have:
\[ r(d+1) = \left\{\begin{array}{ll}
r(d)+1&\mbox{if $r(d)<\binom{d-r(d)}{2}$,}\\
r(d)&\mbox{if $r(d)=\binom{d-r(d)}{2}$.}
\end{array}\right.\]
We claim that $r(d)=\binom{d-r(d)}{2}$ if and only if $d$ is a
triangular number: when $d=\binom{t+1}{2}$ for some $t\geq 0$, we have
\[\binom{t}{2} = \binom{\binom{t+1}{2}-\binom{t}{2}}{2} =  \binom{d-\binom{t}{2}}{2},\]
so $r(d) = \binom{t}{2} = \binom{d-r(d)}{2}$. Conversely,
if $r(d)=\binom{d-r(d)}{2}$, then solving for $d$ we obtain $d =
\binom{d-r(d)+1}{2}$, proving that $d$ is a triangular number.
Therefore, we have:
\[ r(d+1) = \left\{\begin{array}{ll}
r(d) + 1&\mbox{if $d$ is not a triangular number,}\\
r(d) &\mbox{if $d$ is a triangular number.}
\end{array}\right.\]
Since $r(1)=0$, we conclude that $r(d)$ is the number of nontriangular
numbers strictly smaller than~$d$, as claimed. To establish the
formula, note that the value of $r$ at $\binom{t}{2}$ is
$\binom{t-1}{2}$, and therefore $r\left(\binom{t}{2}+s\right) =
\binom{t-1}{2}+(s-1)$ for $0<s<t$, since there are exactly $s-1$ more
nontriangular numbers strictly less than $\binom{t}{2}+s$ than there
are strictly less than $\binom{t}{2}$. And
$\binom{t}{2}+t=\binom{t+1}{2}$, so we also get equality when $s=t$.
\end{proof}

\begin{remark} These alternate descriptions can also be obtained by
  examining sequence  A083920 in~\cite{onlineintseq}; for example, compare the closed formula
  there with $(\ref{eq:maxoverlap})$. In fact, the author first
  realized these descriptions hold by calculating the first few
  values of $r(d)$ directly, and then consulting~\cite{onlineintseq}.
\end{remark}

We can now obtain an upper bound for $\sum\dim(Z_k)$ in terms of
$\dim(X)$, which in turn gives a lower bound for $\dim(X^*)$ in terms
of~$\dim(X)$.

\begin{definition} For $m>0$, we let $f(m)$ denote the largest
  possible value of $\sum\dim(Z_k)$ for a subspace $X$ of $V$ with
  $\dim(X)=m$, for a suitable choice of~$n$. Equivalently,
\[ f(m) = \max\Bigl\{ \dim\bigl(X^n\cap \ker(\Phi)\bigr)\,\Bigm|\, X<V(n),\
  \dim(X)=m\Bigr\}.\]
\end{definition}

\begin{remark} It would appear that this quantity should really be a
  function of $m$ and~$n$, $f(m,n)$. It is easy to verify that
  $f(m,n)\leq f(m,N)$ for any $N\geq n$: if $X$
  is a subspace of $V(n)$, we can also consider it as a subspace of
  $V(N)$ for any $N\geq n$. If the dimension of $X^*$ (with respect to
  $\{\varphi_i\}_{i=1}^n$) is $nm-k$, then 
  the dimension of $X^*$ (with respect to $\{\varphi_i\}_{i=1}^N$) is
  $Nm-k$, so we have
 \[ \dim\bigl(X^n\cap\ker(\Phi_n)\bigr) = \dim\bigl(
  X^N\cap\ker(\Phi_N)\bigr).\] Intuitively, the reason the reverse
inequality also holds is that the largest value of $f(m,n)$ occurs
when the vectors in $Z_i$ use fewer indices rather than more, since
more indices implies a smaller value for $\dim(X\cap\langle
v_{sr}\,|\, i<r<s\leq n\rangle)$, and hence a smaller possible value
for $Z_{j}$ with $j>i$. So the ``best'' strategy for larger
intersection with $\ker(\Phi)$ is to keep $X$ confined to as small a
number of indices as possible.  The proof
below will formalize this intuition, and show that indeed the value of
$f$ depends only on~$m$.
\label{rem:largerndoesnotmatter}
\end{remark}

\begin{theorem} Let $m>0$, and write $m = \binom{T}{2}+s$, $0\leq
  s\leq T$. Then
\[ f(m) =\binom{T}{3} + \binom{s}{2}.\]
\label{th:maxtotaloverlappossible}
\end{theorem}

\begin{remark} Although there is some ambiguity in the expression
  for~$m$, since $\binom{T}{2} + T = \binom{T+1}{2}$, note that the
  values $\binom{T}{3}+\binom{T}{2}$ and $\binom{T+1}{3}+\binom{0}{2}$
  are equal, so the given value of $f(m)$ is well-defined.
\end{remark}

\begin{proof} Assume $s>0$.  First we show that $f(m)\geq
  \binom{T}{3}+\binom{s}{2}$. 

Let $X$ be the $m$-dimensional coordinate subspace of $V(T+1)$ (the
smallest allowable value of $n$) generated by all
$v_{ji}$ with $1\leq i<j\leq T$, and the vectors
$v_{T+1,1},\ldots,v_{T+1,s}$. Then 
\[ X^* = \bigl\langle
\varphi_1(X),\ldots,\varphi_{T}(X),\varphi_{T+1}(X)\bigr\rangle\]
is a coordinate subspace of $W$ generated by the vectors $w_{jik}$,
$1\leq i< j\leq T$, $i\leq k\leq T+1$, and the vectors $w_{T+1,i,k}$
with $1\leq i\leq s$, $i\leq k\leq T+1$. This gives a total of
$2\binom{T}{3} + 3\binom{T}{2} + s(T+1) - \binom{s}{2}$ basis
vectors.  Therefore,
\begin{eqnarray*}
\lefteqn{(T+1)m - \dim(X^*)}\\
 & = & (T+1)\left(\binom{T}{2} + s\right) -
  2\binom{T}{3} - 3\binom{T}{2} - s(T+1) + \binom{s}{2}\\
& = & (T-2)\binom{T}{2} - 2\binom{T}{3} +   \binom{s}{2}\\
& = & \binom{T}{3} + \binom{s}{2},
\end{eqnarray*}
as claimed. As noted in Remark~\ref{rem:largerndoesnotmatter}, this
shows $f(m)\geq \binom{T}{3} + \binom{s}{2}$ for all $n$ which
satisfies $m\leq\binom{n}{2}$.

For the reverse inequality, we will apply induction. Fix 
$m=\binom{T}{2}+s$ with $0<s\leq T$.
 Let $X$ be a subspace of $V(n)$ of dimension~$m$, where $n$ is any
 integer with $m\leq\binom{n}{2}$. We want to show that $\sum
 \dim(Z_i)$ is bounded above by $\binom{T}{3}+\binom{s}{2}$. If all
 $Z_i$ are trivial, this follows. Otherwise, assume $i$ is the
 smallest index with nontrivial $Z_i$, and that $\dim(Z_i)=k>0$. Then
 $k\leq r(m)$, and if $\ell$ is the smallest positive integer such that $k\leq
 \binom{\ell}{2}$ then 
\[ \dim\bigl(X \cap \langle v_{sr}\,|\, i<r<s\leq n\rangle\bigr)\leq
 m-\ell.\]
Thus, $\sum\dim(Z_k) \leq k+ f(m-\ell)$. We want to show
that this expression is bounded above by
$\binom{T}{3}+\binom{s}{2}$. Note that $k\leq
r(m)\leq\binom{m-r(m)}{2}$, so $\ell\leq m-r(m)=T$.

It is easy to show that for $m=1,2,3,4$, and~$5$, all values of the
form $k+f(m-\ell)$, $k\leq r(m)$ and $\ell$ as above are less than or
equal to $\binom{T}{3}+\binom{s}{2}$, which establishes the base case
of the induction.

If $\ell = T = m-r(m)$, then since $k\leq r(m)$ we have
\begin{eqnarray*}
k+f(m-\ell) &\leq &r(m) + f(r(m))\\
 &=& \binom{T-1}{2} + (s-1) +
f\left(\binom{T-1}{2}+(s-1)\right)\\
&=& \binom{T-1}{2} + (s-1) + \binom{T-1}{3} + \binom{s-1}{2}\\
& = & \binom{T}{3} + \binom{s}{2};
\end{eqnarray*}
(this holds even if $s=1$, using induction and as noted above). 
If $\ell<T$, then 
\[\binom{\ell}{2}\leq \binom{T-1}{2} \leq  \binom{T-1}{2}+(s-1) = r(m),\]
and since $k\leq \binom{\ell}{2}$, it is enough to consider the
expression $\binom{\ell}{2}+f(m-\ell)$. When $\ell=T$ this
simplification may not be possible, since we will have
$r(m)<\binom{T}{2}$ whenever $s<T$.  To finish the proof it is
enough to show that for $1<\ell<T$,
\[ \binom{\ell}{2} + f(m-\ell) \leq \binom{T}{3} + \binom{s}{2}.\]
If $2\leq \ell \leq s$, then:
\begin{eqnarray*}
\binom{\ell}{2} + f(m-\ell) & = & \binom{\ell}{2} +
f\left(\binom{T}{2} + (s-\ell)\right)\\
& = & \binom{\ell}{2} + \binom{T}{3} + \binom{s-\ell}{2}\\
&\leq&\binom{T}{3} + \binom{s}{2}.
\end{eqnarray*}
The last inequality follows since $\binom{\ell}{2}+\binom{s-\ell}{2}$
is the number of two element subsets of $\{1,\ldots,s\}$, where either
both elements are less than or equal to $\ell$, or both strictly
larger than $\ell$.

If $s<\ell<T$, then write $\ell = s+a$, $a>0$. We then have
\[ m- \ell = \binom{T}{2} + s - (s+a) = \binom{T-1}{2} + (T-1-a),\]
so
\[ \binom{\ell}{2} + f(m-\ell) = \binom{\ell}{2} + \binom{T-1}{3} +
\binom{T-1-a}{2}.\]
Since $\ell+1-t\leq 0$ and $a>0$, we must have
\[ 6a(s+a+1-T)\leq 0.\]
Rewriting and introducing suitable terms we have:
\[ 6as +3a^2 - 3a - 3T^2 + 9T - 6 + 3T^2 - 9T - 6aT + 9a + 3a^2 +
6\leq 0\]
In turn, this can be rewritten as
\[ 6as + 3a^2 - 3a - 3(T-1)(T-2) + 3(T-a-1)(T-a-2)\leq 0.\]
This gives:
\[ 3(s^2 + 2as + a^2 - s - a) - 3(T-1)(T-2) + 3(T-a-1)(T-a-2)\leq
3(s^2 - s),\]
and so
\[ 3((s+a)^2 - (s+a)) -3(T-1)(T-2) + 3(T-a-1)(T-a-2) \leq 3(s^2-s).\]
Substituting $\ell$ for $s+a$ and adding $T(T-1)(T-2)$ to both sides we have
\[ 3(\ell^2-\ell) + (T-3)(T-2)(T-1) + 3(T-a-1)(T-a-2) \leq T(T-1)(T-2)
+ 3(s^2-s),\]
and dividing through by $6$ yields the desired inequality:
\[\binom{\ell}{2} + \binom{T-1}{3} + \binom{T-1-a}{2} \leq 
\binom{T}{3} + \binom{s}{2}.
\]

We therefore conclude that $f(m)\leq \binom{T}{3}+\binom{s}{2}$, which
completes the proof. Note that indeed, the value of $n$ is not relevant.
\end{proof}

\begin{theorem}
Fix $n>1$ and let $X$ be a subspace of $V$. Write $\dim(X) =
\binom{T}{2}+s$, $0\leq s\leq T$. Then
\[ n\dim(X) - \binom{T}{3} - \binom{s}{2} \leq \dim(X^*) \leq
\min\left\{n\dim(X),\ 2\binom{n+1}{3}\right\}.\]
\label{th:upperandlowerbounds}
\end{theorem}

\begin{proof} The lower bound follows from
$\dim(X^*)\geq n\dim(X) - f(\dim(X))$, and the upper  bound is immediate.
\end{proof}

\begin{corollary} Fix $n>1$ and let $X$ be a subspace of $V$ with
  $\dim(X)=m$. If $\dim(X^*) = nm - k$ and $n+k > f(m+1)$, then $X$ is
  closed.
\end{corollary}

\begin{proof} Suppose $X$ is as in the statement, and let $Y$ be any
  subspace of~$V$ of dimension $m+1$. From the definition
  of $f$ we know that
\[ \dim(Y^*) \geq n(m+1) - f(m+1),\]
so $\dim(Y^*)-\dim(X^*) \geq n+k - f(m+1)>0$.
Therefore every $Y$ strictly larger than $X$ must have
$\dim(X^*)<\dim(Y^*)$, which shows that $X$ is closed by Proposition~\ref{prop:strictlylarger}. 
\end{proof}

\begin{corollary} Fix $n>1$ and let $X$ be a subspace of $V$ with
  $\dim(X)=m$. Write $m = \binom{T}{2} + s$, $0\leq s< T$. If
  $\binom{T}{3}+\binom{s+1}{2}<n$, then $X$ is closed.
\label{cor:Xsmallenough}
\end{corollary}

\begin{proof} This follows from the previous corollary and the formula
  for $f(m+1)$ in Theorem~\ref{th:maxtotaloverlappossible}. 
\end{proof}

For reference, Table~\ref{table:valuesoff} contains the values of
$f(m)$, $3\leq m \leq 50$. Note that $f(1)=f(2)=0$ by Corollary~\ref{cor:maincountingarg}

\begin{table}[t]
\begin{tabular}{||c|c|||c|c|||c|c||}
\hline
$\quad m\quad$&$\quad f(m)\quad$&$\quad m\quad$&$\quad
f(m)\quad$&$\quad m\quad$&$\quad f(m)\quad$\\
\hline\hline
$3$&$1$&$19$&$26$&$35$&$77$\\
\hline
$4$&$1$&$20$&$30$&$36$&$84$\\
\hline
$5$&$2$&$21$&$35$&$37$&$84$\\
\hline
$6$&$4$&$22$&$35$&$38$&$85$\\
\hline
$7$&$4$&$23$&$36$&$39$&$87$ \\
\hline
$8$&$5$&$24$&$38$&$40$&$90$ \\
\hline
$9$&$7$&$25$&$41$&$41$&$94$\\
\hline
$10$&$10$&$26$&$45$&$42$&$99$\\
\hline
$11$&$10$&$27$&$50$&$43$&$105$\\
\hline
$12$&$11$&$28$&$56$&$44$&$112$\\
\hline
$13$&$13$&$29$&$56$&$45$&$120$\\
\hline
$14$&$16$&$30$&$57$&$46$&$120$\\
\hline
$15$&$20$&$31$&$59$&$47$&$121$\\
\hline
$16$&$20$&$32$&$62$&$48$&$123$\\
\hline
$17$&$21$&$33$&$66$&$49$&$126$\\
\hline
$18$&$23$&$34$&$71$&$50$&$130$\\
\hline
\end{tabular}
\caption{Explicit values of $f(m)$, $3\leq m\leq 50$}
\label{table:valuesoff}
\end{table}

Translating back into group theory, we obtain the following:

\begin{theorem} Let $G$ be a group of class at most two and exponent
  $p$, where $p$ is an odd prime. Let ${\rm rank}(G^{\rm ab}) = n$, and let ${\rm
  rank}([G,G])=m$. If $f\bigl(\binom{n}{2}-m+1\bigr)<n$, where $f(k)$ is the
  function in Theorem~\ref{th:maxtotaloverlappossible}, then $G$ is capable.
\label{th:largeenoughcomm}
\end{theorem}

\begin{proof} The subspace $X$ of $V(n)$ corresponding to $G$ has
  dimension $\binom{n}{2}-m$; so the result follows directly from
  Corollary~\ref{cor:Xsmallenough}. 
\end{proof}

\begin{remark} Below we will be able to replace ${\rm rank}(G^{\rm
  ab})$ with ${\rm rank}(G/Z(G))$. 
\end{remark}

\section{Some results.}\label{sec:oldresults}

\begin{lemma} Fix $i,j$, $1\leq i<j\leq n$. If
  $\pi_{ji}(X)=\mathbf{0}$, then $\pi_{ji}(X^{**})=\mathbf{0}$.
\end{lemma}

\begin{proof} Note that $\pi_{jii}(\varphi_k(\mathbf{v}))\neq\mathbf{0}$
 if and only if $k=i$ and $\pi_{ji}(\mathbf{v})\neq \mathbf{0}$. From the
  hypothesis we thus conclude that 
  $\pi_{jii}(X^*)=\mathbf{0}$. This is the same as
  $\pi_{jii}(X^{***})$, so $\pi_{ji}(X^{**})=\mathbf{0}$, as claimed.
\end{proof}

\begin{lemma}[Lemma 5.3 in~\cite{capablep}] If $X$ is a coordinate
  subspace (that is, it is generated by a subset of the $v_{ji}$) then
  $X$ is closed.
\label{lemma:coordinate}
\end{lemma}
\begin{proof} Let $S\subset \{ v_{ji}\,|\, 1\leq i<j\leq n\}$. If
  $X\subset \langle S\rangle$, then by the previous lemma we have
  $X^{**}\subset \langle S \rangle$. If, moreover, $X=\langle
  S\rangle$, then we deduce that $X^{**}\subset X\subset X^{**}$, so
  $X$ is closed.
\end{proof}

We quote the following result without proof:

\begin{lemma}[Lemma 5.7 in~\cite{capablep}]
Let $m$ be fixed, $1<m<n$, and assume that $X_1$ is a
  subspace of $\langle v_{ji}\,|\, 1\leq i<j\leq m\rangle$, and $X_2$
  is a subspace of $\langle v_{ji}\,|\, m<i<j\leq n\rangle$. Then
  $X_1\oplus X_2$ is closed.
\label{lemma:smallplusbig}
\end{lemma}

\begin{lemma} If $\Pi_i(X)=\mathbf{0}$ for some $i$, $1\leq i\leq n$,
  then $X$ is closed.
\end{lemma}

\begin{proof}
By Theorem~\ref{th:descriptionkernelPhi}(iii), if 
$(\mathbf{v}_1,\ldots,\mathbf{v}_n)\in X^n\cap\ker(\Phi)$, then for
$j>i$ we have:
\[\mathbf{0}= \Pi_{i}(\mathbf{v}_j) = \sum_{r=1}^{i-1}\left(-\alpha_{jr}^{(i)}\right)v_{ir}
+ \sum_{r=i+1}^{j-1}\alpha_{jr}^{(i)}v_{ri} + \sum_{r=j+1}^n
\left(-\alpha_{rj}^{(i)}\right)v_{ri},\]
and if $j<i$ then
\[ \mathbf{0} = 
\Pi_{i}(\mathbf{v}_j)  = 
\sum_{r=1}^{j-1}\left(-\alpha_{jr}^{(i)}\right)v_{ir} +
\sum_{r=j+1}^{i-1}\alpha_{rj}^{(i)}v_{ir} +
\sum_{r=i+1}^n\left(-\alpha_{rj}^{(i)}\right)v_{ri}.\]
Thus we conclude that $\alpha_{rs}^{(i)}=0$ for all $1\leq r<s\leq n$.
Therefore, $\mathbf{v}_i=\mathbf{0}$. In particular,
the intersection of $\varphi_i(X)$ and $\langle \varphi_j(X)\,|\,
j\neq i\rangle$ is trivial, so $\varphi_i^{-1}(X^*)=X$. Therefore,
$X^{**}\subset \varphi_i^{-1}(X^*) = X$, proving $X$ is
closed.
\end{proof}

In~\cite{capablep}*{Lemma~5.9} we proved a special case of the following result:

\begin{lemma}
Let $n>1$ be an integer. Suppose  $I$ is a proper nonempty subset of
$\{1,\ldots,n\}$, and let $J= \{1,\ldots,n\}-I$. Let 
\begin{eqnarray*}
V_I &=& \langle v_{ji}\,|\, 1\leq i<j\leq n;\ i,j\in I\}\\
V_J &=& \langle v_{ji}\,|\, 1\leq i<j\leq n;\ i,j\in J\}\\
V_{(I,J)} & = & \langle v_{ji}\,|\,1\leq i<j\leq n;\ \mbox{exactly one
  of $i,j$ is in $I$, one in $J$}\rangle.
\end{eqnarray*}
Let $X_I$ be a subspace of $V_I$, $X_J$ be a subspace of
$V_J$, and let
\[ X = X_I \oplus X_J \oplus V_{(I,J)}.\]
Then 
\[ X^{**} = \left(\bigcap_{i\in
  I}\varphi_i^{-1}\Bigl(\bigl\langle \varphi_k(X_I)\,\bigm|\,k\in I\bigr\rangle\Bigr)\right)
    \oplus \left(\bigcap_{j\in
  J}\varphi_j^{-1}\Bigl(\bigl\langle \varphi_k(X_J)\,\bigm|\, k \in
  J\bigr\rangle\Bigr)\right)\oplus V_{(I,J)}.\]
In particular, $X$ is
  $\{\varphi_k\}_{k=1}^n$-closed if and only if $X_I$ is
  $\{\varphi_i\}_{i\in I}$-closed and $X_J$ is $\{\varphi_j\}_{j\in
  J}$-closed.
\end{lemma}

\begin{proof} Let $W'$ be the subspace of $W$ generated by all $w_{jik}$ in which
exactly one or two of $i,j,k$ are in~$I$. 
First, we claim that $W'\subset X^*$.
Suppose that exactly one of $i,j,k$ lies in~$I$. If $j\in I$, then
$v_{ji} \in V_{(I,J)}$, so $\varphi_k(v_{ji})=w_{jik}\in X^*$. If
$i\in I$, then again $v_{ji}\in V_{(I,J)}$, so again we conclude that
$w_{jik}\in X^*$. Finally, if $k\in I$, then either $v_{jk}$ or $v_{kj}$ lies
in $V_{(I,J)}$ (whichever is appropriate); the image under $\varphi_i$
of this vector is $\pm (w_{jik}-w_{kij})$. Since $w_{kij}\in
X^*$ by the argument above it follows
that $w_{jik}\in X^*$ as desired.
The case when exactly two if $i,j,k$ lie in $I$ follows by symmetry,
since exactly one of them will lie in~$J$.

As a second step, we claim that
\[ X^* = \Bigl\langle \varphi_i(X_I)\,\Bigm|\, i\in I\Bigr\rangle
\oplus \Bigl\langle \varphi_j(X_J)\,\Bigm|\, j\in J\Bigr\rangle \oplus
W'.\]
Indeed, the right hand side is contained in $X^*$.
It suffices to show that $\varphi_j(X_I)\subset W'$ for every $j\in J$
(the symmetrical argument will show it for $\varphi_i(X_J)$, $i\in I$).
Since $\pi_{abc}\left(\varphi_j(X_I)\right)\neq \{\mathbf{0}\}$ 
only if exactly two of $a,b,c$ are in~$I$ and either $b$ or $c$
are equal to~$j$, this establishes our second claim.

Thirdly, if we let $W_I$ (resp.~$W_J$) denote the subspace of~$W$
generated by all $w_{jik}$ in which all of $i,j,k$ lie in~$I$
(resp.~$i,j,k$ lie in~$J$), then we claim that
\[ X^*\cap W_I = \langle \varphi_i(X_I)\,|\,i\in I\rangle,\quad\mbox{and}\quad
X^*\cap W_J = \langle \varphi_j(X_J)\,|\, j\in J\rangle.\]
Clearly, $\varphi_i(X_I)\subset X^*\cap W_I$ for all $i\in I$. For the
converse inclusion simply note that $X^* = X_I^* + X_J^* +
V_{(I,J)}^*$, and both $X_J^*$ and $V_{(I,J)}^*$ are disjoint from
$W_I$. As $\varphi_k(X_I)$ is contained in $W'$ when $k\in J$ and
contained in $W_I$ when $k\in I$, the claim now follows. A symmetric
argument holds for $X^*\cap W_J$.

Finally, we establish the equality in the lemma: let
\[ \mathbf{v} \in \left(\bigcap_{i\in
  I}\varphi_i^{-1}\Bigl(\bigl\langle \varphi_k(X_I)\,\bigm|\,k\in I\bigr\rangle\Bigr)\right)
     \oplus \left(\bigcap_{j\in
  J}\varphi_j^{-1}\Bigl(\bigl\langle \varphi_k(X_J)\,\bigm|\, k \in
  J\bigr\rangle\Bigr)\right)\oplus V_{(I,J)}.\]
Then we can write $\mathbf{v} = \mathbf{x}_I \oplus \mathbf{x}_{(I,J)}
  \oplus \mathbf{x}_J$, where
\begin{eqnarray*}
\mathbf{x}_I & \in & \bigcap_{i\in I}\varphi_i^{-1}\Bigl(\bigl\langle
\varphi_k(X_I)\,|\,k\in I\bigr\rangle\Bigr),\\
\mathbf{x}_J & \in & \bigcap_{j\in J}\varphi_j^{-1}\Bigl(\bigl\langle
\varphi_k(X_J)\,|\, k\in J\bigr\rangle\Bigr),\\
\mathbf{x}_{(I,J)} & \in & V_{(I,J)}.
\end{eqnarray*}
Since $V_{(I,J)}\subset X \subset X^{**}$, to show that $\mathbf{v}\in
X^{**}$ it is enough to show that $\mathbf{x}_I$ and $\mathbf{x}_J$
are both in~$X^{**}$. And indeed: notice that
$\pi_{ab}(\mathbf{x}_I)\neq\mathbf{0}$ only if $a,b\in I$. Therefore,
$\varphi_j(\mathbf{x}_I)\in W'\subset X^*$ for every $j\in J$. And by
construction we have $\varphi_i(\mathbf{x}_I)\in X^*$ for every $i\in
I$. Thus $\varphi_k(\mathbf{x}_I)\in X^*$ for all~$k$. Symmetrically,
$\varphi_k(\mathbf{x}_J)\in X^*$ for all~$k$ as well. Therefore, each
of $\mathbf{x}_I$ and $\mathbf{x}_J$ lie in $X^{**}$, as desired.

For the converse inclusion, let $\mathbf{x}\in X^{**}$. We can write
$\mathbf{x}=\mathbf{v}_I \oplus\mathbf{v}_J\oplus \mathbf{v}_{(I,J)}$,
with $\mathbf{v}_I\in V_I$, $\mathbf{v}_J\in V_J$, and
$\mathbf{v}_{(I,J)}\in V_{(I,J)}$.  We must have
$\varphi_k(\mathbf{v})\in X^*$ for every $k$. If $k\in I$, then 
we have 
\[\varphi_k(\mathbf{v}_I) \in X^*\cap W_I = \Bigl\langle \varphi_i(X_I)\,\Bigm|\,
i\in I\Bigr\rangle,\qquad\mbox{and}\qquad \varphi_k(\mathbf{v}_J),\varphi_k(\mathbf{v}_{(I,J)})\in W'.\]
And if $k\in J$, then
\[ \varphi_k(\mathbf{v}_J)\in X^*\cap W_J = \Bigl\langle \varphi_j(X_J)\,|\Bigm|\,
j\in J\Bigr\rangle,\qquad\mbox{and}\qquad
\varphi_k(\mathbf{v}_I),\varphi_k(\mathbf{v}_{(I,J)})\in W'.\]
Therefore, 
\begin{eqnarray*}
\mathbf{v}_I & \in & \bigcap_{i\in I}\varphi_i^{-1}\Bigl(\bigl\langle
\varphi_k(X_I)\,|\,k\in I\bigr\rangle\Bigr),\\
\mathbf{v}_J & \in & \bigcap_{j\in J}\varphi_j^{-1}\Bigl(\bigl\langle
\varphi_k(X_J)\,|\, k\in J\bigr\rangle\Bigr),
\end{eqnarray*}
as desired. This proves the equality.

The final clause of the lemma follows from the observation that
$X$ will be closed if and only if 
\[X_I = \bigcap_{i\in I} \varphi_i^{-1}\Bigl(\bigl\langle
\varphi_k(X_I)\,\bigm|\, k\in I\bigr\rangle\Bigr)\ \mbox{and}\ X_J =
\bigcap_{j\in J}\varphi_j^{-1}\Bigl(\bigl\langle
\varphi_k(X_J)\,\bigm|\, k\in J\bigr\rangle\Bigr).\]
\end{proof}

Notice that when we interpret the subspace $V_{(I,J)}$ as the subgroup
of commutators we are making trivial, this subgroup
``says'' that each $x_i$ commutes with each $x_j$ in~$G$, where $i\in I$ and
$j\in J$. Translating back into group theoretic terms we obtain the
following:

\begin{corollary} Let $G_1$ and $G_2$ be two noncyclic $p$-groups of
  class at most two and exponent~$p$. Then $G=G_1\oplus G_2$ is
  capable if and only if each of $G_1$ and~$G_2$ are capable.
\end{corollary}

The corollary is clearly not true if we drop the ``noncyclic''
hypothesis, since a cyclic group of order~$p$ is not capable, but the
direct sum of two cyclic groups of order~$p$ is capable. To understand
why we must make the distinction, note that when $|I|=1$, $V_I$ is
trivial so $X_I$ is perforce $\{\varphi_i\}_{i\in I}$-closed by vacuity; but our
developement always assumes $n>1$. The symmetric situation happens
when $|I|=n-1$. However, what we obtain in this instance is the
following result, that is perhaps one of the most unexpected results
from~\cite{capablep}: 

\begin{theorem}[Lemma~5.9 in~\cite{capablep}] Let $n>2$, let $X$ be a
  subspace of $V$ such that $\Pi_n(X)=\mathbf{0}$,
  and let $X' = X\oplus \langle
  v_{ni}\,|\, 1\leq i\leq n-1\rangle$. Then
  $X'$ is $\{\varphi_i\}_{i=1}^{n}$-closed if and only if $X$ is
  $\{\varphi_i\}_{i=1}^{n-1}$-closed.
\label{th:cancelcentral}
\end{theorem}

Note that the subspace $\langle v_{ni}\,|\,1\leq i\leq n-1\rangle$ is
none other than $\langle u_n\rangle^*$ (with respect to
$\{\psi_i\}_{i=1}^n$). The choice of $u_n$ is merely one of
convenience; if we invoke symmetry, we obtain:

\begin{corollary} Let $n>2$ and let $X'$ be a subspace of $V(n)$. If
  there exists a vector $\mathbf{u}\in U(n)$, $\mathbf{u}\neq\mathbf{0}$, such
  that $\langle \mathbf{u}\rangle^*$ is contained in $X'$, then there
  exists a subspace $X$ of $V(n-1)$ such that  
  $X'$ is $\{\varphi_i\}_{i=1}^n$ closed if and only if $X$ is
  $\{\varphi_i\}_{i=1}^{n-1}$ closed, and $\dim(X)=\dim(X')-n+1$.
\label{cor:cancelanycentral}
\end{corollary}

The usefulness of these theorems becomes apparent when we translate it
back into group theoretic terms:

\begin{corollary}  Let $K$ be a finite noncyclic group of exponent $p$ and
  class $2$, and let $G=K\oplus C_p$, where $C_p$ is the
  cyclic group of order~$p$. Then $G$ is capable if and only if $K$ is
  capable.
\end{corollary}

Since every finite group of class at most two and exponent $p$ may be
written as $K\oplus C_p^r$ for some $r\geq 0$ and $K$ satisfying
$[K,K]=Z(K)$, we conclude that:

\begin{corollary} Let $G$ be a $p$-group of class at most two and
  exponent~$p$. Write $G=K\oplus C_p^r$, where $K$ satisfies
  $[K,K]=Z(K)$. Then $G$ is capable if and only if (i) $K$ is nontrivial
  and capable, or (ii) $K$ is trivial and $r\geq 2$.
\end{corollary}

This allows us to ignore the ``extra'' elements in the center of~$G$
that do not come from commutators. Noting that $G/Z(G)\cong K^{\rm
  ab}$ and $[G,G]=[K,K]$, we can, for example, 
strengthen Theorem~\ref{th:largeenoughcomm}
by replacing the rank of
$G^{\rm ab}$ by the rank of $G/Z(G)$. We then have the following:

\begin{corollary}
Let $G$ be a $p$-group of class at most two and exponent~$p$. Let ${\rm
  rank}(G/Z(G))=n$, and ${\rm rank}([G,G])=m$.
If $f\bigl(\binom{n}{2} - m + 1\bigr) < n$, where $f(k)$ is the function in
  Theorem~\ref{th:maxtotaloverlappossible}, then 
then $G$ is capable.
\label{cor:suffcardcond}
\end{corollary}

In~\cite{heinnikolova}, the authors established a minimum size for the
commutator subgroup of a capable group~$G$ of class two and
exponent~$p$ which satisfies $[G,G]=Z(G)$; namely, they proved:

\begin{theorem}[Theorem~1 in~\cite{heinnikolova}] Suppose that $G$ is a group of
  exponent~$p$ which satisfies $[G,G]=Z(G)$. If $G$ is capable and
  $[G,G]$ is of rank~$m$, then $G/Z(G)$ is of rank at most $2m +
  \binom{m}{2}$.
\label{th:heinnik}
\end{theorem}

Our development above now shows that if $G$ is not cyclic, then
Theorem~\ref{th:heinnik} remains true if we weaken the hypothesis that
$[G,G]=Z(G)$ to $[G,G]\subset Z(G)$ (i.e., class at most two). Note also that
Theorem~\ref{th:heinnik} gives a \textit{necessary} condition. Combining it with
Corollary~\ref{cor:suffcardcond}, we have:

\begin{theorem} Let $G$ be a noncylic group of exponent~$p$ and class
  at most two. Let ${\rm rank}(G/Z(G))=n$ and ${\rm rank}([G,G])=m$.
A necessary condition for the capablity of~$G$ is that $n$ and $m$ satisfy
$n\leq 2m+\binom{m}{2}$.
A sufficient condition for the capability of~$G$ is that $n$ and $m$ satisfy
$f\left(\binom{n}{2}-m+1\right) < n$,
where $f(k)$ is the function in Theorem~\ref{th:maxtotaloverlappossible}.
\end{theorem}

The theorem can be interpreted as saying that a group of class exactly
$2$ and exponent $p$ is capable only if it is ``nonabelian enough'':
the necessary condition shows that $m$ cannot be too small relative to
$n$ (so there must be ``enough'' nontrivial commutators), while the
sufficient condition shows that if the commutator subgroup is large
enough then the group will necessarily be capable.

\section{Central elements and the kernel of $\Phi$.}\label{sec:centralelements}

Let $F$ be a field, and let $k,n$ be integers. The Grassmannian
$Gr(k,n)$ is the set of all $k$-dimensional subspaces of $F^n$. This
set is in fact an algebraic variety with very rich structure; for
example, $Gr(1,n)$ is $\mathbb{P}^{n-1}$, projective $(n-1)$-space over $F$.

Corollary~\ref{cor:cancelanycentral} hints at the importance of
subspaces of~$V$ the form $Z^*$, where
$Z$ is a subspace of~$U$, for the question of
capability. In this section we will explore some of the connections
between these subspaces and $\ker(\Phi)$. 

\begin{lemma}
Fix $n>1$ and let $Z$ be a subspace of~$U$. If $\dim(Z)=k$, then
$\dim(Z^*) = k(n-k)+\binom{k}{2} = kn - \binom{k+1}{2}$.
\label{lemma:sizeofZstar}
\end{lemma}

\begin{proof} Let $\mathbf{z}_1,\ldots,\mathbf{z}_k$ be a basis for $Z$. Let
  $\mathbf{u}_{k+1},\ldots,\mathbf{u}_{n}$ be vectors of~$U$ that
  extend this to a basis for $U$. Then $Z^*\cong Z\wedge U$, and a
  basis for $Z\wedge U$ is given by all elements of the form
  $\mathbf{z}_j\wedge\mathbf{z}_i$, with $1\leq i<j\leq k$, and all
  elements of the form $\mathbf{z}_j\wedge \mathbf{u}_r$, with $1\leq
  j\leq k$ and $k+1\leq r\leq n$.
\end{proof}

\begin{lemma} Fix $n>1$, and let $k$ be a positive integer. If
  $k+1<n$, then the map $Gr(k,U)\to Gr(nk-\binom{k+1}{2},V)$ given by
  $Z\mapsto Z^*$ is one-to-one.
\label{lemma:onetoone}
\end{lemma}

\begin{proof} If $Z^*=Z'^*$, then any subspace $Y\subset \langle
  Z,Z'\rangle$ satisfies $Y^*\subset Z^*$. Thus, it is enough to show
  that if $\dim(Y)=k+1$, then $\dim(Y^*)>\dim(Z^*)$. This is
  equivalent to verifying the inequality $n(k+1)-\binom{k+2}{2}>nk -
  \binom{k+1}{2}$, which holds exactly when $k+1<n$. 
\end{proof}

The special case $k=1$ is of
particular interest:

\begin{definition} We define $\Psi\colon\mathbb{P}^{n-1}=Gr(1,U)\to Gr(n-1,V)$
  as the map that sends $[\alpha_1 : \cdots : \alpha_n]$ to
  $(\alpha_1u_1+\cdots+\alpha_nu_n)\wedge U$; i.e., the subspace
  generated by
\[
\mathbf{v}_1  = \sum_{j=1}^n \alpha_j v_{j1},\quad \mathbf{v}_2 =
\sum_{j=1}^n \alpha_j v_{j2}, \quad \ldots\;,\quad \mathbf{v}_n =
\sum_{j=1}^n \alpha_jv_{jn},\]
with $v_{ij}=-v_{ji}$ and $v_{ii}=\mathbf{0}$.
\label{def:defofPsi}
\end{definition}

The subspaces which are images under this map are
closely connected to $\ker(\Phi)$. We explore this connection in the
next series of results.

\begin{lemma} Let $\mathbf{p}=[\alpha_1 :\cdots :
    \alpha_n]\in\mathbb{P}^{n-1}$; if $\alpha_i\neq 0$, then a basis
    for $\Psi(\mathbf{p})$ is given by
    $\mathbf{v}_1,\ldots,\widehat{\mathbf{v}_i},\ldots,\mathbf{v}_n$,
    i.e., taking all $\mathbf{v}_j$ and omitting the vector
    $\mathbf{v}_i$.
\label{lemma:forcommsalwaysnminusone}
\end{lemma}

\begin{proof} By Lemma~\ref{lemma:sizeofZstar} it is enough to show
    that the given list is linearly independent.
Indeed, if $\beta_1\mathbf{v}_1 + \cdots +
\widehat{\beta_i\mathbf{v}_i} + \cdots + \beta_n\mathbf{v}_n =
\mathbf{0}$ (where as usual we use $\widehat{\beta_i\mathbf{v}_i}$ to
mean we are omitting that summand), then for fixed $j\neq i$ we have:
\[ \pi_{ji}(\beta_1\mathbf{v}_1 + \cdots +
\widehat{\beta_i\mathbf{v}_i} + \cdots + \beta_n\mathbf{v}_n)  =
\beta_j\alpha_i v_{ij}=\mathbf{0},\]
and since $\alpha_i\neq 0$ we conclude that $\beta_j=0$.
\end{proof}

\begin{lemma} Fix $n>1$, let $\mathbf{p}\in\mathbb{P}^{n-1}$,
  $\mathbf{p}=[\alpha_1:\cdots:\alpha_n]$, and let
  $\mathbf{v}\in\Psi(\mathbf{p})$. If $\alpha_i\neq 0$ and
  $\Pi_i(\mathbf{v})=\mathbf{0}$, then $\mathbf{v}=\mathbf{0}$.
\label{lemma:ithprojofPsi}
\end{lemma}

\begin{proof} Write $\mathbf{v}=\beta_1\mathbf{v}_1 + \cdots +
  \widehat{\beta_i\mathbf{v}_i} + \cdots +
  \beta_n\mathbf{v}_n$. If $j\neq i$, then 
$\pi_{ji}(\mathbf{v}) = \beta_j\alpha_i v_{ij}$. Since $\alpha_i\neq 0$
  by hypothesis, it follows that $\beta_j=0$. Therefore,
  $\mathbf{v}=\mathbf{0}$, as claimed.
\end{proof}

\begin{lemma} Fix $n>1$ and let $\mathbf{p}\in\mathbb{P}^{n-1}$. Then
  $\Psi(\mathbf{p})^n\cap\ker(\Phi)$ is
  trivial.
\label{lemma:PsiintersectPhi}
\end{lemma}

\begin{proof} Let $(\mathbf{w}_1,\ldots,\mathbf{w}_n)\in
  \ker(\Phi)\cap(\Psi(\mathbf{p}))^n$. Write
$\mathbf{p}=[\alpha_1:\cdots:\alpha_n]$.
By Theorem~\ref{th:descriptionkernelPhi}(ii),
$\Pi_j(\mathbf{w}_j)=\mathbf{0}$ for all~$j$. If $\alpha_i\neq 0$, then by
Lemma~\ref{lemma:ithprojofPsi} we must have
$\mathbf{w}_i=\mathbf{0}$. If we then let $j\neq i$, then by
Theorem~\ref{th:descriptionkernelPhi}(iii) we also have
$\Pi_i(\mathbf{w}_j)=\mathbf{0}$, hence $\mathbf{w}_j=\mathbf{0}$ for
all $j$, as claimed.
\end{proof}

\begin{theorem} Fix $n>1$, and let $\mathbf{p}\in
  \mathbb{P}^{n-1}$, $\mathbf{v}\in V$, and $X = \langle
  \Psi(\mathbf{p}),\mathbf{v}\rangle$. Then $X^n\cap\ker(\Phi)$ is
  nontrivial if and only if $\mathbf{v}\notin \Psi(\mathbf{p})$.
\label{th:ifcentralplusonethenoverlap}
\end{theorem}

\begin{proof}  Let $\mathbf{p}=[\alpha_1 : \cdots :\alpha_n]$ be
  an element of $\mathbb{P}^{n-1}$. Let
  $\mathbf{v}_1,\ldots,\mathbf{v}_n$ be the vectors
\[ \mathbf{v}_i = \sum_{j=1}^n \alpha_{j}v_{ji};\]
we know that if $\alpha_i\neq 0$, then
$\mathbf{v}_1,\ldots,\widehat{\mathbf{v}_i},\ldots,\mathbf{v}_n$ is a
basis for $\Psi(\mathbf{p})$. 

The ``only if'' clause follows from
Lemma~\ref{lemma:PsiintersectPhi}. For the ``if'' clause, 
let $\mathbf{v}$ be a vector of $V$ which is not in
$\Psi(\mathbf{v})$, with
\[ \mathbf{v}=\sum_{1\leq i<j\leq n}
\!\!\!\!\alpha_{ji}v_{ji}.\]
Let $\alpha_{ji}=-\alpha_{ij}$ and
$\alpha_{ii}=0$ for simplicity.
We claim that
\[ \mathbf{w} = \varphi_1\left(\alpha_1\mathbf{v} +
\sum_{j=1}^n\alpha_{j1}\mathbf{v}_j\right) + \cdots +
\varphi_n\left(\alpha_n\mathbf{v} +
\sum_{j=1}^{n}\alpha_{jn}\mathbf{v}_n\right)=\mathbf{0}.\]
In general we have:
\begin{eqnarray*}
\pi_{m\ell}\left(\alpha_u\mathbf{v} + \sum_{j=1}^n
  \alpha_{ju}\mathbf{v}_j\right)
& = & \pi_{m\ell}\left( \sum_{1\leq i<j\leq n}\!\!\!\!\!\!\alpha_u\alpha_{ji}v_{ji} +
  \sum_{j=1}^n\sum_{i=1}^n\alpha_{ju}\alpha_iv_{ij}\right)\\
& = & \alpha_u\alpha_{m\ell}v_{m\ell} +
  (\alpha_{mu}\alpha_{\ell}v_{\ell m} + \alpha_{\ell
  u}\alpha_mv_{m\ell})\\
& = & \left(\alpha_{m\ell}\alpha_u + 
  \alpha_{um}\alpha_{\ell} + \alpha_{\ell u}\alpha_m\right)v_{m\ell}.
\end{eqnarray*}
Let $r,s,t$ be integers, $1\leq r<s\leq n$,
$s\leq t\leq n$. We want to prove that
$\pi_{srt}(\mathbf{w})=\mathbf{0}$. The $(srt)$ coefficient of
$\mathbf{w}$ is equal to the $(tr)$ coefficient of
$\alpha_t\mathbf{v}+\sum\alpha_{jt}\mathbf{v}_j$ when $s=t$; to the
sum of the $(sr)$ coefficient of
$\alpha_t\mathbf{v}+\sum\alpha_{jt}\mathbf{v}_j$ and the $(st)$
coefficient of $\alpha_r\mathbf{v}+\sum\alpha_{jr}\mathbf{v}_j$ when
$s>t$; and to the differentce of the $(sr)$ coefficient of
$\alpha_t\mathbf{v}+\sum\alpha_{jt}\mathbf{v}_j$ and the $(ts)$
coefficient of $\alpha_r\mathbf{v}+\sum\alpha_{jr}\mathbf{v}_j$ when
$s<t$. 

Therefore, if $s=t$ then we have:
$\pi_{trt}(\mathbf{w})  = \left(\alpha_{tr}\alpha_t +
\alpha_{rt}\alpha_t + \alpha_{tt}\alpha_r\right)w_{trt} = \mathbf{0}$,
(since $\alpha_{tt}=0$ and $\alpha_{rt}=-\alpha_{tr}$). If $s>t$ then
we have: 
\[\pi_{srt}(\mathbf{w}) =\left( (\alpha_{sr}\alpha_t +
\alpha_{rt}\alpha_s + \alpha_{ts}\alpha_r) + (\alpha_{st}\alpha_r +
\alpha_{tr}\alpha_s + \alpha_{rs}\alpha_t)\right)w_{srt} =
\mathbf{0}.\]
And finally if $s<t$ then we have:
\[\pi_{srt}(\mathbf{w}) = \left( (\alpha_{sr}\alpha_t +
\alpha_{rt}\alpha_s + \alpha_{ts}\alpha_r) - (\alpha_{ts}\alpha_r +
\alpha_{sr}\alpha_t + \alpha_{rt}\alpha_s)\right) w_{srt} =
\mathbf{0}.\]
Thus, $\mathbf{w}=\mathbf{0}$ and so
\[ \left(
\alpha_1\mathbf{v}+\sum_{j=1}^n\alpha_{j1}\mathbf{v}_j,\;\ldots\;,\alpha_n\mathbf{v}+\sum_{j=1}^n\alpha_{jn}\mathbf{v}_j\right)\]
lies in $X^n\cap\ker(\Phi)$.

We claim that this element is nonzero if and only if $\mathbf{v}\notin
\Psi(\mathbf{p})$. Indeed, if the element is trivial, then since
$\alpha_i\neq 0$ for some $i$, the $i$-th component will yield an
expression for $\mathbf{v}$ as a linear combination of
$\mathbf{v}_1,\ldots,\mathbf{v}_n$,
proving that $\mathbf{v}$ lies in $\Psi(\mathbf{p})$. Conversely,
assume that $\mathbf{v}\in\Psi(\mathbf{p})$. Assuming $\alpha_i\neq
0$, then we can write
\[ \mathbf{v} = \sum_{j=1}^n \beta_j\mathbf{v}_j =
\sum_{j=1}^n\sum_{k=1}^n \beta_j\alpha_k v_{kj},\quad\mbox{where
  $\beta_i=0$}.\]
Then for any pair $a,b$, $1\leq a<b\leq n$, we have $\alpha_{ba} =
  \beta_{a}\alpha_b - \beta_b\alpha_a$, and therefore:
\begin{eqnarray*}
\lefteqn{\pi_{m\ell}\left(\alpha_u\mathbf{v} + \sum_{j=1}^n
  \alpha_{ju}\mathbf{v}_j\right)
 =  \left(\alpha_{m\ell}\alpha_u + \alpha_{\ell u}\alpha_m +
  \alpha_{um}\alpha_{\ell}\right)v_{m\ell}}\\
& = & \left( \beta_m\alpha_{\ell}\alpha_u -
  \beta_{\ell}\alpha_m\alpha_u + \beta_{\ell}\alpha_u\alpha_m -
  \beta_u\alpha_{\ell}\alpha_m + \beta_u\alpha_m\alpha_{\ell} -
  \beta_m\alpha_u\alpha_{\ell}\right)v_{m\ell}\\
& = & \mathbf{0},
\end{eqnarray*}
as claimed. In particular, if $\mathbf{v}\notin \Psi(\mathbf{p})$,
then $X^n\cap \ker(\Phi)$ is nontrivial, proving the theorem.
\end{proof}

\section{The case $n=4$.}\label{sec:nequalfour}

In this section we will settle the 
  case of $n$-generator groups of class two and exponent~$p$, with
  $n\leq 4$.

From Table~\ref{table:valuesoff} and
  Corollary~\ref{cor:suffcardcond} we deduce that for $n\leq 4$,
  all subspaces $X$ of $V(n)$ are closed, with the possible
  exception of $n=4$ and $\dim(X)=5$ (in that case, we have
  $f(\dim(X)+1)=  4$, so the corollary does not apply; if
  $\dim(X)=6$ and $n=4$, then $X=V$ which is trivially closed).

Thus, we may restrict our attention to the case $n=4$ and
$\dim(X)=5$. The only two possibilities are $X^{**}=X$ and
$X^{**}=V$. Since $X^{**}=V$ if and only if $X^*=W$, we have:

\begin{prop}
Let $n=4$, and let $X$ be a $5$-dimensional subspace of~$V$. Then $X$
is closed if and only if $X^*$ is a proper subspace
of~$W$. Equivalently, $X$ is closed if and only if $X^n\cap\ker(\Phi)$
is nontrivial.
\end{prop}

\begin{proof} For the last assertion, note that $n\dim(X) =
  20=\dim(W)$, so $X^*\neq W$ if and only if
  $\dim(X^n\cap\ker(\Phi))>0$.
\end{proof}

By Proposition~\ref{prop:kerPhifornfour}, if
$\mathbf{v}=(\mathbf{v}_1,\mathbf{v}_2,\mathbf{v}_3,\mathbf{v}_4)\in
\ker(\Phi)$, then either
$\mathbf{v}=(\mathbf{0},\mathbf{0},\mathbf{0},\mathbf{0})$, or else
the subspace spanned by
$\mathbf{v}_1,\mathbf{v}_2,\mathbf{v}_3,\mathbf{v}_4$ has
dimension~$3$. Since $\ker(\Phi)$ is of dimension $\binom{n}{3}$, in
this case dimension~$4$, the
one-dimensional subspaces of $\ker(\Phi)$ correspond to points in $\mathbb{P}^3$.
Thus we obtain a map from $\mathbb{P}^3$ to $Gr(3,V)$. Explicitly:

\begin{definition} Let
  $\mathbf{q}=[\alpha_{123}:\alpha_{124}:\alpha_{134}:\alpha_{234}]\in\mathbb{P}^3$.
  Then $\Upsilon(\mathbf{q})$ is defined to be the element of
  $Gr(3,V)$ spanned by
\begin{eqnarray*}
\mathbf{v}_1 & = & \alpha_{123} v_{32} + \alpha_{124} v_{42} + \alpha_{134} v_{43},\\
\mathbf{v}_2 & = & -\alpha_{123}v_{31} - \alpha_{124} v_{41} + \alpha_{234} v_{43},\\
\mathbf{v}_3 & = & \alpha_{123} v_{21} - \alpha_{134} v_{41} - \alpha_{234} v_{42},\\
\mathbf{v}_4 & = & \alpha_{124} v_{21} + \alpha_{134} v_{31} + \alpha_{234} v_{32}.
\end{eqnarray*}
\end{definition}

Note that $(\mathbf{v}_1,\mathbf{v}_2,\mathbf{v}_3,\mathbf{v}_4)\in
\ker(\Phi)$, and each nonzero element of $\ker(\Phi)$ corresponds to a
point in $\mathbb{P}^3$.  We are using triples of integers as indices
because we are ``really'' working over $\mathbb{P}^{\binom{4}{3}-1}$.

Thus we have:

\begin{corollary} Let $n=4$ and let $X$ be a $5$-dimensional subspace
  of~$V$. Then $X$ is closed if and only if $X$ contains
  $\Upsilon(\mathbf{q})$ for some $\mathbf{q}\in \mathbb{P}^3$.
\end{corollary}

Using this notation, we can rephrase
Lemma~\ref{lemma:PsiintersectPhi} and
Theorem~\ref{th:ifcentralplusonethenoverlap} for $n=4$ as follows: 

\begin{lemma} Fix $n=4$ and let $\mathbf{p}\in\mathbb{P}^{3}$. Then
  $\Psi(\mathbf{p})\neq \Upsilon(\mathbf{q})$ for any $\mathbf{q}\in \mathbb{P}^3$.
\label{lemma:PsiintersectPhinfour}
\end{lemma}

\begin{proof} From Lemma~\ref{lemma:PsiintersectPhi} we know that
  $\Psi(\mathbf{p})$ cannot contain $\Upsilon(\mathbf{q})$ for any
  $\mathbf{q}$; since $\Upsilon(\mathbf{q})$ is of dimension
  $3$, the statement follows.
\end{proof}

\begin{lemma} 
Let $n=4$, and let $\mathbf{p}\in
  \mathbb{P}^{3}$, $\mathbf{v}\in V$, and $X = \langle
  \Psi(\mathbf{p}),\mathbf{v}\rangle$. Then there exists
  $\mathbf{q}\in\mathbb{P}^3$ such that $\Upsilon(\mathbf{q})$ is
  contained in $X$ if and only if $\mathbf{v}\notin \Psi(\mathbf{p})$.
\end{lemma}

In the case $n=4$ we can prove the symmetric result:

\begin{theorem} Let $n=4$ and let $\mathbf{q}\in \mathbb{P}^3$,
  $\mathbf{v}\in V$, and $X = \langle
  \Upsilon(\mathbf{q}),\mathbf{v}\rangle$. 
Then 
  there exists $\mathbf{p}\in\mathbb{P}^3$ such that
  $\Psi(\mathbf{p})$ is contained in $X$ if and only if
  $\mathbf{v}\notin\Upsilon(\mathbf{q})$. 
\label{th:conversenfour}
\end{theorem}

\begin{proof} Let
  $\mathbf{q}=[\alpha_{123} : \alpha_{124} : \alpha_{134} : \alpha_{234}]\in\mathbb{P}^3$.
Then $\Upsilon(\mathbf{q})$ is generated by the vectors:
\begin{eqnarray*}
\mathbf{v}_1 & = & \alpha_{123} v_{32} + \alpha_{124} v_{42} + \alpha_{134} v_{43},\\
\mathbf{v}_2 & = & -\alpha_{123}v_{31} - \alpha_{124} v_{41} + \alpha_{234} v_{43},\\
\mathbf{v}_3 & = & \alpha_{123} v_{21} - \alpha_{134} v_{41} - \alpha_{234} v_{42},\\
\mathbf{v}_4 & = & \alpha_{124} v_{21} + \alpha_{134} v_{31} + \alpha_{234} v_{32}.
\end{eqnarray*}

The ``only if'' clause follows from
Lemma~\ref{lemma:PsiintersectPhinfour}.  For the ``if'' clause,
let $\mathbf{v}\in V$ be the vector given by
\[ \mathbf{v} = \!\!\!\!\sum_{1\leq i<j\leq n}\!\!\!\! \alpha_{ji}v_{ji},\]
and let $\mathbf{p}\in\mathbb{P}^3$ be given by $\mathbf{p}=[\beta_1 :
  \beta_2 : \beta_3 : \beta_4]$, where
\begin{eqnarray*}
\beta_1 & = & \alpha_{123}\alpha_{41} - \alpha_{124}\alpha_{31} +
\alpha_{134}\alpha_{21}\\
\beta_2 & = & \alpha_{123}\alpha_{42} - \alpha_{124}\alpha_{32} +
\alpha_{234}\alpha_{21}\\
\beta_3 & = & \alpha_{123}\alpha_{43} - \alpha_{134}\alpha_{32} +
\alpha_{234}\alpha_{31}\\
\beta_4 & = & \alpha_{124}\alpha_{43} - \alpha_{134}\alpha_{42} +
  \alpha_{234}\alpha_{41}.
\end{eqnarray*}
(Technically, we need to show that not all $\beta_i$ are equal to $0$
to justify that $\mathbf{p}\in\mathbb{P}^3$; we do this below).

Let $\mathbf{w}_1$, $\mathbf{w}_2$, $\mathbf{w}_3$, and
$\mathbf{w}_4$ be the four vectors as in
Definition~\ref{def:defofPsi}, i.e., 
\[ \mathbf{w}_i  = \sum_{j=1}^4 \beta_j v_{ji},\quad i=1,2,3,4,\]
where as usual we let $v_{ij}=-v_{ji}$ and $v_{ii}=\mathbf{0}$.
It is straightforward to verify that:
\begin{eqnarray*}
\mathbf{w}_1 & = & \alpha_{234}\mathbf{v} - \alpha_{43}\mathbf{v}_2 +
\alpha_{42}\mathbf{v}_3 - \alpha_{32}\mathbf{v}_4,\\
\mathbf{w}_2 & = & -\alpha_{134}\mathbf{v} + \alpha_{43}\mathbf{v}_1 -
\alpha_{41}\mathbf{v}_3 + \alpha_{31}\mathbf{v}_4,\\
\mathbf{w}_3 & = & \alpha_{124}\mathbf{v} - \alpha_{42}\mathbf{v}_1 +
\alpha_{41}\mathbf{v}_2 - \alpha_{21}\mathbf{v}_4,\\
\mathbf{w}_4 & = & -\alpha_{123}\mathbf{v} + \alpha_{32}\mathbf{v}_1 -
\alpha_{31}\mathbf{v}_2 + \alpha_{21}\mathbf{v}_3.
\end{eqnarray*}
The pattern in the above is as follows: for $\mathbf{w}_i$, we take
$\alpha_{abc}\mathbf{v}$, where $a<b<c$ and all three are different
from~$i$, plus the sum of all three cyclic permutations of the indices:
$\alpha_{ab}\mathbf{v}_c + \alpha_{bc}\mathbf{v}_a +
\alpha_{ca}\mathbf{v}_b$,
with the usual proviso that $\alpha_{ji}=-\alpha_{ij}$. In the case of
$\mathbf{w}_2$ and $\mathbf{w}_4$ we further multiply everything by
$-1$.

We claim that all $\beta_i$ are equal to $0$ if and only if
$\mathbf{v}\in\Upsilon(\mathbf{q})$. 
Indeed: if $\beta_i=0$ for each $i$, then $\mathbf{w}_j=\mathbf{0}$
for each $j$. At least one $\alpha_{abc}$ is nonzero, so if $d\notin
\{a,b,c\}$, $1\leq d\leq 4$, then the formula for $\mathbf{w}_d$
expresses $\mathbf{v}$ as a linear combination of three of the
$\mathbf{v}_j$; hence $\mathbf{v}\in \Upsilon(\mathbf{q})$. 

Conversely, if $\mathbf{v}\in\Upsilon(\mathbf{q})$, we want to show
that all $\beta_i$ are equal to $0$. Write:
\[ \mathbf{v} = \gamma_1\mathbf{v}_1 + \gamma_2\mathbf{v}_2 +
\gamma_3\mathbf{v}_3 + \gamma_4\mathbf{v}_4,\]
where we fix $a,b,c$, $1\leq a<b<c\leq 4$, such that $\alpha_{abc}\neq
0$, and set $\gamma_d=0$, where $d$ is the element of $\{1,2,3,4\}$
not in $\{a,b,c\}$. With this proviso, we have:
\begin{eqnarray*}
\alpha_{21} & = & \gamma_3\alpha_{123} + \gamma_4\alpha_{124},\\
\alpha_{31} & = & -\gamma_2\alpha_{123} + \gamma_4\alpha_{134},\\
\alpha_{41} & = & -\gamma_2\alpha_{124} - \gamma_3\alpha_{134},\\
\alpha_{32} & = & \gamma_1\alpha_{123} + \gamma_4\alpha_{234},\\
\alpha_{42} & = & \gamma_1\alpha_{124} - \gamma_3\alpha_{234},\\
\alpha_{43} & = & \gamma_1\alpha_{134} + \gamma_2\alpha_{234}.
\end{eqnarray*}
From this, $\beta_1=\beta_2=\beta_3=\beta_4=0$ readily follows. 

Therefore, $\mathbf{p}$ is a well-defined point in $\mathbb{P}^3$ if
and only if $\mathbf{v}\notin\Upsilon(\mathbf{q})$, from which the
theorem follows.
\end{proof}

So we obtain:

\begin{corollary} Let $n=4$, and let $X$ be a $5$-dimensional subspace
  of~$V(4)$. Then $X$ is closed if and only if there exists
  $\mathbf{p}\in\mathbb{P}^3$ such that $\Psi(\mathbf{p})\subset X$.
\end{corollary}

Translated into group theory clarifies the situation:

\begin{theorem} Let $G$ be a group of class two and exponent~$p$, $p$
  an odd prime, and assume that $G^{\rm ab}$ is of rank~$4$. If
  $[G,G]$ is of rank~$1$, then $G$ is capable if and only if
  $Z(G)/[G,G]$ is nontrivial; that is, if and only if $G$ is not
  extra-special.
\end{theorem}

Recall that we say a group is $k$-generated if it can be generated by
$k$ elements, though it may need fewer. We obtain:

\begin{theorem} Let $G$ be a $4$-generated group of class at most $2$
  and exponent an odd prime~$p$. Then $G$ is one and only of:
\begin{itemize}
\item[(i)] Cyclic and nontrivial;
\item[(ii)] Extra special of order $p^5$ and exponent~$p$;
\item[(iii)] Capable.
\end{itemize}
\end{theorem}

\begin{remark} Here is an alternative proof of
  Theorem~\ref{th:conversenfour} which is entirely geometrical, due to
  David McKinnon~\cite{persdave}. 
The
  maps $\Psi\colon \mathbb{P}^3\to Gr(3,V)$ and
  $\Upsilon\colon\mathbb{P}^3\to Gr(3,V)$ are both regular maps; that
  is, they are defined everywhere, and are locally (in the Zariski
  topology) determined by
  rational functions on the coordinates; in fact, in this case, by
  linear functions on the coordinates. We define two subsets of the
  variety $Gr(4,V)\times \mathbb{P}^3$ by
\begin{eqnarray*}
A & = & \Bigl\{ (X,\mathbf{p})\,\Bigm|\, \Psi(\mathbf{p})\subset X
\Bigr\},\\
B & = & \Bigl\{ (X,\mathbf{q})\,\Bigm|\, \Upsilon(\mathbf{q})\subset X
\Bigr\}.
\end{eqnarray*}
Since both $\Psi$ and $\Upsilon$ are regular, it follows that both $A$
and $B$ are closed subvarieties of $Gr(4,V)\times \mathbb{P}^3$. The
projections 
\begin{align*}
p_1 &\colon Gr(4,V)\times \mathbb{P}^3\to Gr(4,V),\\
p_2 &\colon Gr(4,V)\times \mathbb{P}^3\to \mathbb{P}^3,
\end{align*}
induce maps from each of $A$ and $B$ into $Gr(4,V)$ and
$\mathbb{P}^3$. The maps to $\mathbb{P}^3$ are surjections, and the
fibers all have dimension~$2$ because the fiber over $\mathbf{p}$
(resp.~$\mathbf{q}$) is 
the set of all $4$-dimensions subspaces of $V$ which contain the
$3$-dimensional subspace $\Psi(\mathbf{p})$
(resp.~$\Upsilon(\mathbf{p})$). This set is isomorphic to the set
of lines in the quotient space $V/\Psi(\mathbf{p})$
(resp.~$V/\Upsilon(\mathbf{p})$), which is in turn isomorphic to
$\mathbb{P}^2$, so it is $2$-dimensional.

The maps are also smooth, so we have smooth maps of fiber dimension
$2$ over a smooth $3$-dimensional variety, hence $A$ and $B$ are both
of dimension $3+2=5$.

Consider now the projections to $Gr(4,6)$; 
we know that $p_1(A)$ and $p_1(B)$ are irreducible subvarieties of
$Gr(4,6)$ of dimension at most~$5$, and that $p_1(A)$ is contained in
$p_1(B)$ (since by Theorem~\ref{th:ifcentralplusonethenoverlap}, if
$(X,\mathbf{p})$ is in $A$, then there exists $\mathbf{q}$ such that
$(X,\mathbf{q})\in B$). If we can show that $p_1(A)$ has dimension
exactly $5$, then the irreducibility of $p_1(B)$ will imply that
$p_1(B)=p_1(A)$, which will show that if $X$ is any $4$-dimensional
subspace of~$V$, and $(X,\mathbf{q})\in B$, then there exists
$\mathbf{p}$ such that $(X,\mathbf{p})\in A$, which is what
Theorem~\ref{th:conversenfour} states.

To show that $p_1(A)$ has dimension exactly $5$, it is enough to show
that it is generically finite; for this it is enough to show that
there is at least one $X\in Gr(4,6)$ such that $p_1^{-1}(X)$ is
nonempty and finite. But in fact we know that $p_1^{-1}(X)$ has at
most one element, since $\mathbf{p}\neq\mathbf{q}$ implies that
$\langle\Psi(\mathbf{p}),\Psi(\mathbf{q})\rangle$ has dimension $5$ by
Lemma~\ref{lemma:sizeofZstar}. Thus, Theorem~\ref{th:conversenfour}
follows.
\end{remark}

\begin{remark} Unfortunately, The analogue of
  Theorem~\ref{th:conversenfour} does not hold for $n=3$ or for
  $n>4$. At least for odd~$n$, we can consider the point $\mathbf{q}$
  which has components equal to $0$ except for those associated to the
  triples $(1,2,3)$, $(1,4,5)$, $(1,6,7)$, etc. Then it is easy to
  verify that $\Upsilon(\mathbf{q})$ is of dimension~$n$, generated by
  all pairs $v_{j1}$ with $j>1$, and the vector 
$\mathbf{v}=v_{32} + v_{54} + v_{76} + \cdots + v_{2n+1,2n}$. Thus,
the subspace $\Upsilon(\mathbf{q})$ contains $\langle u_1\rangle^*$. 
\end{remark}

\section*{Acknowledgements}
I thank David McKinnon for many stimulating conversations, and
for his help in finding the formula for $f(m)$.
Part of this work was conducted while I was on a brief visit to
the Department of Pure Mathematics of the University of Waterloo, at
the invitation of Prof.~McKinnon; 
I am very grateful to him for the invitation, and to the University of
Waterloo for the hospitality and facilities I received
there.

\section*{References}
\bibliographystyle{amsalpha}

\end{document}